\documentclass[10pt]{amsart} 
\usepackage{amssymb}
 
\newtheorem{thm1}{Theorem}[section]

\newtheorem{rem1}[thm1]{Remark}

\newtheorem{cor1}[thm1]{Corollary}

\newtheorem{prop1}[thm1]{Proposition}
\newtheorem{ex1}[thm1]{Example}
\newtheorem{not1}[thm1]{Notation}

\begin{document}

\title[]
{On complete intersection projective closures of monomial curves}
\author[]{Anargyros Katsabekis}
\address {Department of Mathematics, University of Ioannina, 45110 Ioannina, Greece} \email{katsampekis@uoi.gr}

\begin{abstract} We study the complete intersection property of projective closures of affine monomial curves. For monomial curves in the four-dimensional affine space, we give a complete characterization of when the defining ideals of their projective closures are complete intersections. We also investigate when the projective closures are arithmetically Cohen--Macaulay.
\end{abstract}
\maketitle

\section{Introduction}

Toric varieties are fundamental objects in algebraic geometry and commutative algebra. Among their key properties, the complete intersection property is particularly important, as it implies that the associated coordinate ring is Cohen--Macaulay. A toric variety is called a {\em complete intersection} if its associated toric ideal is a {\em complete intersection}, that is, if it can be generated by a number of polynomials equal to its height.

Determining when the projective closure of a toric variety is a complete intersection is a complicated problem, because homogenization can change both the defining equations and the number of minimal generators of the toric ideal. While the complete intersection property has been thoroughly studied for affine monomial curves (see, e.g., \cite{Bresinsky75}, \cite{Her}), much less is known about their projective closures.

The paper \cite{BG2} studies the complete intersection property for projective closures of certain classes of monomial curves, while \cite{BG} considers it for simplicial projective toric varieties. However, until now, it has been an open problem to determine when the projective closure of a monomial curve in the four-dimensional affine space is a complete intersection. Note that the projective closure can be a complete intersection only if the affine monomial curve itself is. Indeed, dehomogenizing a generating set of the defining ideal of the projective closure yields a generating set of the affine toric ideal.

We classify all complete intersection projective closures of monomial curves in the four-dimensional affine space. Our approach combines a detailed analysis of the minimal binomial generators of the defining toric ideals with Gr\"obner basis techniques.

Let ${\bf a}=(a_{1},\ldots,a_{n})$ be a sequence of distinct positive integers with ${\rm gcd}(a_{1},\ldots,a_{n})=1$, and let $C({\bf a})$ be the affine monomial curve parametrized by $x_{i}=t^{a_i}$ for $1 \leq i \leq n$. Consider the polynomial ring $K[x_{1},\ldots,x_{n}]$ over any field $K$. The {\em toric ideal} of $C({\bf a})$, denoted by $I({\bf a})$, is the kernel of the $K$-algebra homomorphism $\phi \colon K[x_{1},\ldots,x_{n}] \to K[t]$ given by $\phi(x_{i})=t^{a_i}$ for all $1 \leq i \leq n$. The ideal $I({\bf a})$ is generated by binomials of the form $x_{1}^{u_1} \cdots x_{n}^{u_n}-x_{1}^{v_1} \cdots x_{n}^{v_n}$ such that $u_{1}a_{1}+\cdots+u_{n}a_{n}=v_{1}a_{1}+\cdots+v_{n}a_{n}$ (see, for instance, \cite[Lemma 4.1]{Sturmfels95}). By \cite[Lemma 4.2]{Sturmfels95}, the height of $I({\bf a})$ is $n-1$, and hence
$I({\bf a})$ is a complete intersection if and only if it is minimally generated by $n-1$ binomials.

Set $d={\rm max}\{a_1,\ldots,a_n\}$. The {\em homogenization} $I^{h}({\bf a})$ of $I({\bf a})$ with respect to the variable $x_{0}$ is the kernel of the $K$-algebra homomorphism $\psi \colon K[x_{0},x_1,\ldots,x_n] \to K[s,t]$, defined by $\psi(x_{0})=s^{d}$ and $\psi(x_{i})=s^{d-a_{i}}t^{a_i}$ for all $1 \leq i \leq n$. Then $I^{h}({\bf a})$ is a complete intersection if and only if it is minimally generated by $n-1$ binomials. 

The {\em projective closure} $\overline{C({\bf a})} \subset \mathbb{P}^{n}$ of $C({\bf a})$ is the projective monomial curve parametrized by $[s:t] \mapsto [s^{d}:s^{d-a_{1}}t^{a_1}: \cdots: s^{d-a_{n}}t^{a_n}]$ and its defining homogeneous ideal is $I^{h}({\bf a})$. The curve $\overline{C({\bf a})}$ is called {\em arithmetically Cohen--Macaulay} if $I^{h}({\bf a})$ is Cohen--Macaulay. This property has been extensively studied; see, for example, \cite{HerSta}, \cite{Katsa1}, \cite{PatilRoberts2003}, and \cite{SSS}.

In Section 2, we observe that when $n=4$ and $I({\bf a})$ is a complete intersection toric ideal, its minimal generating set falls into one of three distinct cases (see Theorem \ref{CasesMinimal}). For the first case, we establish necessary and sufficient conditions for the projective closure to be arithmetically Cohen--Macaulay (see Theorems \ref{ACM1}, \ref{ACM2}, \ref{ACM3}, and \ref{ACM4}) and for $I^{h}({\bf a})$ to be a complete intersection (see Theorems \ref{BAS1} and \ref{CI23}, and Corollaries \ref{CorCase11CI} and \ref{CorCase1CI}).

In Section 3, we consider the second case and determine necessary and sufficient conditions under which $I^{h}({\bf a})$ is a complete intersection (see Theorems \ref{Veryfirst}, \ref{Basic4}, and \ref{Comp2}, and Corollaries \ref{Cor12} and \ref{Cor13}). We also establish necessary and sufficient conditions for $\overline{C({\bf a})}$ to be arithmetically Cohen--Macaulay in certain subcases (see Theorems \ref{ACM5} and \ref{ACM6}). Moreover, Proposition~\ref{ACM7} gives an infinite family of monomial curves for which the projective closure is arithmetically Cohen--Macaulay, while the corresponding homogeneous toric ideal is not a complete intersection and can have arbitrarily many minimal generators.

In Section 4, we study the third case and derive explicit criteria for $I^{h}({\bf a})$ to be a complete intersection (see Theorems \ref{Basic5}, \ref{Basic23}, \ref{basicvfff}, \ref{Basic24}, \ref{Basic25}, and \ref{Basic26}).

\section{Projective closures in case 1}

In this section, we study in detail the first case arising from the structure of the minimal generators. For any $f \in K[x_1,\ldots,x_n]$, we denote by $f^{h}$ its homogenization with respect to $x_0$. The following remark will be used repeatedly.

\begin{rem1}\label{VeryBasicRemark} {\rm (1) Let $I\subset K[x_1,\ldots,x_n]$ be an ideal generated by $f_1,\ldots,f_r$. If the initial monomials of the generators
$f_1,\ldots,f_r$
are pairwise relatively prime, then $\{f_1,\ldots,f_r\}$ is a Gr\"obner basis for $I$ with respect to $<$ by Buchberger's criterion \cite[Proposition~4, Chapter~2]{CLO}.\\ (2) If $\{f_1,\ldots,f_r\}$ is a Gr\"obner basis for an ideal $I$ with respect to a graded monomial order $<$, then $\{f_1^h,\ldots,f_r^h\}$ is a Gr\"obner basis for $I^h$ with respect to the graded monomial order $<_0$ extending $<$ and satisfying $x_0<_0x_i$ for all $i\ge1$ (see \cite[Chapter~8]{CLO}).} \end{rem1}

Let ${\bf a}=(a_1,\ldots,a_4)$ be a sequence of distinct relatively prime positive integers. 
\begin{not1} {\rm Throughout the paper, let $c_i$ denote the smallest positive integer such that $c_{i}a_{i} \in \sum_{j \neq i} \mathbb{N}a_{j}$, for each $i=1,\ldots,4$, where $\mathbb{N}$ is the set of nonnegative integers.}
\end{not1}

\begin{thm1} \label{CasesMinimal} Let $I({\bf a})$ be a complete intersection toric ideal, where ${\bf a}=(a_1,\ldots,a_4)$ is a sequence of distinct positive integers with ${\rm gcd}(a_1,\ldots,a_4)=1$. Then after permuting the variables, if necessary, there exists a minimal system of binomial generators $S$ of $I({\bf a})$ of the following form \begin{enumerate} \item[{\rm CASE 1:}] If $c_{1}a_{1}=c_{2}a_{2}$ and $c_{3}a_{3}=c_{4}a_{4}$, then $S=\{x_{2}^{c_2}-x_{1}^{c_1}, x_{3}^{c_3}-x_{4}^{c_4}, x_{1}^{d_1}x_{2}^{d_2}-x_{3}^{d_3}x_{4}^{d_4}\}$. \item[{\rm CASE 2:}] If $c_{1}a_{1}=c_{2}a_{2}=c_{3}a_{3} \neq c_{4}a_{4}$, then $S=\{x_{2}^{c_2}-x_{1}^{c_1}, x_{3}^{c_3}-x_{1}^{c_1}, x_{4}^{c_4}-x_{1}^{d_1}x_{2}^{d_2}x_{3}^{d_3}\}$, where at least two of $d_1$, $d_2$, and $d_3$ are nonzero. \item[{\rm CASE 3:}] If $c_{1}a_{1}=c_{2}a_{2}$ and $c_{i}a_{i} \neq c_{j}a_{j}$ for all $\{i,j\} \neq \{1,2\}$, then $S=\{x_{2}^{c_2}-x_{1}^{c_1}, x_{3}^{c_3}-x_{1}^{b_1}x_{2}^{b_2}, x_{4}^{c_4}-x_{1}^{d_1}x_{2}^{d_2}x_{3}^{d_3}\}$, where both $b_1$ and $b_2$ are nonzero, and at least two of $d_1$, $d_2$, and $d_3$ are nonzero.
\end{enumerate}
\end{thm1}

{\em \noindent Proof.} We begin by noting that if $I({\bf a})$ is a complete intersection, then one of the following holds: either $c_1a_1=c_2a_2$ and $c_3a_3=c_4a_4$; or $c_{1}a_{1}=c_{2}a_{2}=c_{3}a_{3} \neq c_{4}a_{4}$; or $c_{1}a_{1}=c_{2}a_{2}$ and $c_{i}a_{i} \neq c_{j}a_{j}$ for all $\{i,j\} \neq \{1,2\}$ (see \cite{KO}). By \cite[Proposition 3.2]{Shi}, a minimal generating set of $I({\bf a})$ is of the form $\{x_{2}^{c_2}-x_{1}^{c_1}, x_{3}^{c_3}-x_{4}^{c_4}, x_{1}^{d_1}x_{2}^{d_2}-x_{3}^{d_3}x_{4}^{d_4}\}$ in Case 1, and $\{x_{2}^{c_2}-x_{1}^{c_1}, x_{3}^{c_3}-x_{1}^{b_1}x_{2}^{b_2}, x_{4}^{c_4}-x_{1}^{d_1}x_{2}^{d_2}x_{3}^{d_3}\}$ in Cases 2 and 3. In Case 2, since $c_{1}a_1=c_{2}a_2=c_{3}a_{3}$, we can set $b_{1}=c_1$ and $b_2=0$, so that the second generator is $x_{3}^{c_3}-x_{1}^{c_1}$. We may also assume that at least two of $d_1$, $d_2$, $d_3$ are nonzero; otherwise, if, for instance, $d_{1}=d_{2}=0$, then $d_{3}>c_3$, and using $x_{3}^{c_3}-x_{1}^{c_1}$, we can replace $x_{4}^{c_4}-x_{3}^{d_3}$ with $x_{4}^{c_4}-x_{1}^{c_1}x_{3}^{d_3-c_3} \in I({\bf a})$. In Case 3, we may assume that both $b_1$ and $b_2$ are nonzero; indeed, if, for instance, $b_2=0$, then $b_1>c_1$, and using $x_2^{c_2}-x_1^{c_1}$ we can replace $x_3^{c_3} - x_1^{b_1}x_2^{b_2}$ with $x_3^{c_3}-x_1^{b_1-c_1}x_2^{b_2+c_2} \in I({\bf a})$. Similarly, we may assume that at least two of $d_1$, $d_2$, $d_3$ are nonzero. \hfill $\square$ \par\vspace{0.5\baselineskip}

 For the rest of this section, let $I({\bf a})$ be a complete intersection toric ideal minimally generated by $x_{2}^{c_2}-x_{1}^{c_1}$, $x_{3}^{c_3}-x_{4}^{c_4}$, and $x_{1}^{d_1}x_{2}^{d_2}-x_{3}^{d_3}x_{4}^{d_4}$. {\em Without loss of generality, we fix the assumptions $a_{4}={\rm max}\{a_1,\ldots,a_4\}$, $c_2>c_1$, $d_2 \leq c_2$, and $d_3 \leq c_3$.} Indeed, if $d_2>c_2$, we may write $d_{2}=kc_2+l$ for some integers $k \geq 1$ and $0 \leq l<c_2$, and replace $x_{1}^{d_1}x_{2}^{d_2}-x_{3}^{d_3}x_{4}^{d_4}$ with $x_{1}^{kc_1+d_1}x_{2}^{l}-x_{3}^{d_3}x_{4}^{d_4} \in I({\bf a})$. Similarly, we may assume $d_3 \leq c_3$. Finally, since $a_4>a_3$, we have $c_4<c_3$.

Throughout this section, let $<$ be the degree reverse lexicographic order on $K[x_1, \ldots, x_4]$ with $x_1 > x_2 > x_3 > x_4$. We begin by considering a subcase in which the ideal $I^{h}({\bf a})$ is a complete intersection. \vspace{-0.25\baselineskip}

\begin{thm1} \label{BAS1} If $d_2=c_2$ and $d_3=c_3$, then $I^{h}({\bf a})$ is a complete intersection.
\end{thm1}
\vspace{-0.25\baselineskip}
{\em \noindent Proof.} The ideal $I({\bf a})$ is also generated by the binomials $x_2^{c_2}-x_1^{c_1}$, $x_3^{c_3}-x_4^{c_4}$, and $x_1^{c_1+d_1}-x_4^{c_4+d_4}$.
Their initial monomials are
$x_2^{c_2}$,
$x_3^{c_3}$,
and
$x_1^{c_1+d_1}$,
which are pairwise relatively prime.
Hence, by Remark~\ref{VeryBasicRemark}(1), these binomials form a
Gr\"obner basis for $I({\bf a})$ with respect to $<$.
By Remark~\ref{VeryBasicRemark}(2), their homogenizations form a
Gr\"obner basis for $I^h({\bf a})$ with respect to $<_0$.
Hence $I^h({\bf a})$ is a complete intersection.
\hfill $\square$ \par\vspace{0.5\baselineskip}

The following example illustrates the subcase considered in Theorem~\ref{BAS1}.

\begin{ex1} {\rm Let ${\bf a}=(132, 60, 165, 220)$. Then $I({\bf a})$ is generated by $x_{2}^{11}-x_{1}^{5}$, $x_3^{4}-x_{4}^{3}$, and $x_{2}^{11}-x_{3}^{4}$. Here, $d_{2}=c_2=11$ and $d_{3}=c_3=4$. By Theorem \ref{BAS1}, the ideal $I^{h}({\bf a})$ is a complete intersection.}
\end{ex1}

Next, we consider another subcase in which $I^{h}({\bf a})$ is a complete intersection.

\begin{thm1} \label{CI23} Suppose that $d_2=0$ and $d_{3}=c_{3}$. Then $I^{h}({\bf a})$ is a complete intersection.
\end{thm1}
{\em \noindent Proof.} Let
$f_1=x_2^{c_2}-x_1^{c_1}$,
$f_2=x_3^{c_3}-x_4^{c_4}$,
and
$g=x_1^{d_1}-x_4^{c_4+d_4}\in I({\bf a})$. Since $x_1^{d_1}-x_3^{c_3}x_4^{d_4}=g-x_4^{d_4}f_2$, the ideal $I({\bf a})$ is generated by $f_1$, $f_2$, and $g$.
The initial monomials of $f_1$, $f_2$, and $g$ are
$x_2^{c_2}$,
$x_3^{c_3}$,
and
$x_1^{d_1}$,
respectively, and are pairwise relatively prime.
Hence, by Remark~\ref{VeryBasicRemark}(1),
$\{f_1,f_2,g\}$ is a Gr\"obner basis for $I({\bf a})$ with respect to $<$.
By Remark~\ref{VeryBasicRemark}(2),
$\{f_1^{h},f_2^{h},g^{h}\}$ is a Gr\"obner basis for $I^{h}({\bf a})$
with respect to $<_0$.
Thus $I^{h}({\bf a})$ is a complete intersection.
\hfill $\square$ \par\vspace{0.5\baselineskip}

We now obtain several necessary and sufficient conditions under which the arithmetically Cohen--Macaulay property holds.

\begin{thm1} \label{ACM1} Suppose that $d_2=c_2$ and $d_3<c_3$. Then $\overline{C({\bf a})}$ is arithmetically Cohen--Macaulay if and only if $c_{1}+d_{1} \geq d_{3}+d_{4}$.
\end{thm1}

{\em \noindent Proof.} ($\Leftarrow$) Let $f_{1}=x_2^{c_2}-x_1^{c_1}$, $f_{2}=x_3^{c_3}-x_4^{c_4}$, and $f_{3}=x_1^{c_1+d_{1}}-x_3^{d_3}x_4^{d_4}\in I({\bf a})$. Since $x_{1}^{d_1}x_{2}^{c_2}-x_{3}^{d_3}x_{4}^{d_4}=x_{1}^{d_1}f_1+f_3$, the ideal $I({\bf a})$ is generated by the binomials $f_1$, $f_2$, and $f_3$. Then ${\rm in}_{<}(f_{1})=x_2^{c_2}$, ${\rm in}_{<}(f_{2})=x_3^{c_3}$, and ${\rm in}_{<}(f_{3})=x_1^{c_1+d_1}$ are pairwise relatively prime. Hence, by Remark~\ref{VeryBasicRemark}(1), $\{f_1,f_2,f_3\}$ is a Gr\"obner basis for $I({\bf a})$ with respect to $<$. Since none of these initial monomials is divisible by $x_4$, it follows from \cite[Theorem~2.2]{HerSta} that $\overline{C({\bf a})}$ is arithmetically Cohen--Macaulay.\\
($\Rightarrow$) Assume that $c_{1}+d_{1}<d_{3}+d_{4}$. Since $\overline{C({\bf a})}$ is arithmetically Cohen--Macaulay and ${\rm in}_{<}(f_{3})=x_3^{d_{3}}x_4^{d_4}$, by \cite[Theorem 2.2]{HerSta} there exists a binomial $B \in I({\bf a})$ such that ${\rm in}_{<}(B)$ divides $x_3^{d_3}$; in particular ${\rm in}_{<}(B)=x_{3}^{e}$ for some $e \leq d_3$. Since $d_3<c_3$, this gives $e<c_3$, contradicting the minimality of $c_3$. \hfill $\square$ \par\vspace{0.5\baselineskip}

As an application of Theorem \ref{ACM1}, we obtain the following corollary. 

\begin{cor1} \label{CorCase11CI} Suppose that $d_2=c_2$ and $d_{3}<c_{3}$. Then the ideal $I^{h}({\bf a})$ is a complete intersection if and only if $c_{1}+d_{1} \geq d_{3}+d_{4}$.
\end{cor1}

{\em \noindent Proof.} ($\Leftarrow$) By the proof of Theorem~\ref{ACM1}, $\{f_i \mid 1\le i\le 3\}$ is a Gr\"obner basis for $I({\bf a})$ with respect to $<$. Hence, by Remark~\ref{VeryBasicRemark}(2), $\{f_i^{h} \mid 1\le i\le 3\}$ is a Gr\"obner basis for $I^{h}({\bf a})$ with respect to $<_0$. Thus $I^{h}({\bf a})$ is a complete intersection.\\ ($\Rightarrow$) If $I^{h}({\bf a})$ is a complete intersection, then $\overline{C({\bf a})}$ is arithmetically Cohen--Macaulay. By Theorem \ref{ACM1}, it follows that $c_1+d_{1} \geq d_{3}+d_{4}$. \hfill $\square$ \par\vspace{0.5\baselineskip}

\begin{thm1} \label{ACM2} Suppose that $0<d_2<c_2$ and $d_{3}=c_{3}$. Then $\overline{C({\bf a})}$ is arithmetically Cohen--Macaulay if and only if $c_{1}+d_{1} \geq c_{2}-d_{2}+c_{4}+d_{4}$.
\end{thm1}

{\em \noindent Proof.} ($\Leftarrow$) Consider the set $T=\{f_1=x_2^{c_2}-x_1^{c_1},\; f_2=x_3^{c_3}-x_4^{c_4},\; f_3=x_1^{d_1}x_2^{d_2}-x_4^{c_4+d_4},\; f_4=x_1^{c_1+d_1}-x_2^{c_2-d_2}x_4^{c_4+d_4}\} \subseteq I({\bf a})$. Since $c_2>c_1$ and $c_1+d_1 \geq c_2-d_2+c_4+d_4$, the initial monomials with respect to $<$ are ${\rm in}_{<}(f_{1})=x_2^{c_2}$, ${\rm in}_{<}(f_{2})=x_3^{c_3}$, ${\rm in}_{<}(f_{3})=x_1^{d_1}x_2^{d_2}$, and ${\rm in}_{<}(f_4)=x_1^{c_1+d_{1}}$. Since ${\rm in}_{<}(f_1)$ and ${\rm in}_{<}(f_2)$ are relatively prime, it follows that $S(f_{1},f_{2}) \stackrel{T}{\longrightarrow} 0$. Similarly, $S(f_{1},f_4) \stackrel{T}{\longrightarrow} 0$, $S(f_{2},f_{3}) \stackrel{T}{\longrightarrow} 0$, and $S(f_{2},f_4) \stackrel{T}{\longrightarrow} 0$. We have that $S(f_{1},f_{3})=-f_4 \stackrel{f_4}{\longrightarrow} 0$ and $S(f_{3},f_4)=x_4^{c_4+d_4}(x_{2}^{c_2}-x_{1}^{c_1}) \stackrel{f_1}{\longrightarrow} 0$. Thus $T$ is a Gr\"obner basis for $I({\bf a})$ with respect to $<$. Since none of the initial monomials is divisible by $x_4$, it follows from \cite[Theorem~2.2]{HerSta} that $\overline{C({\bf a})}$ is arithmetically Cohen--Macaulay.\\
($\Rightarrow$) Assume that $c_{1}+d_{1}<c_{2}-d_{2}+c_{4}+d_{4}$. Since $\overline{C({\bf a})}$ is arithmetically Cohen--Macaulay and ${\rm in}_{<}(f_4)=x_2^{c_2-d_2}x_4^{c_4+d_4}$, by \cite[Theorem~2.2]{HerSta} there exists a binomial $B \in I({\bf a})$ such that ${\rm in}_{<}(B)$ divides $x_2^{c_{2}-d_2}$, contradicting the minimality of $c_2$. \hfill $\square$

\begin{thm1} \label{ACM3} Suppose that $0<d_2<c_2$ and $d_{3}<c_{3}$. Then $\overline{C({\bf a})}$ is arithmetically Cohen--Macaulay if and only if $c_{1}+d_{1} \geq c_{2}-d_{2}+d_{3}+d_{4}$.
\end{thm1}

{\em \noindent Proof.} ($\Leftarrow$) The proof proceeds via Gr\"obner basis techniques, as in the proof of Theorem \ref{ACM2}. Since $c_2 > c_1$ and $c_1 + d_1 \geq c_2 - d_2 + d_3 + d_4$, it follows that $d_1 + d_2 > d_3 + d_4$. Let $f_1=x_2^{c_2}-x_1^{c_1}, f_2=x_3^{c_3}-x_4^{c_4}, f_3=x_1^{d_1}x_2^{d_2}-x_3^{d_3}x_4^{d_4}, f_4=x_1^{c_1+d_{1}}-x_2^{c_{2}-d_2}x_3^{d_3}x_4^{d_4} \in I({\bf a})$. With respect to $<$, the initial monomials of the above binomials are $x_2^{c_2}$, $x_3^{c_3}$, $x_1^{d_1}x_2^{d_2}$, and $x_1^{c_1+d_{1}}$, respectively. As in the proof of Theorem~\ref{ACM2}, all $S$-polynomials reduce to zero, so $\{f_i \mid 1 \leq i \leq 4\}$ is a Gr\"obner basis of $I({\bf a})$ with respect to $<$. Since $x_4$ does not divide any initial monomial, $\overline{C({\bf a})}$ is arithmetically Cohen--Macaulay by \cite[Theorem 2.2]{HerSta}.\\
($\Rightarrow$) Assume that $c_{1}+d_{1}<c_{2}-d_{2}+d_{3}+d_{4}$. Since $\overline{C({\bf a})}$ is arithmetically Cohen--Macaulay and ${\rm in}_{<}(f_4)=x_2^{c_{2}-d_2}x_3^{d_3}x_4^{d_4}$, by \cite[Theorem 2.2]{HerSta} there exists a binomial $B \in I({\bf a})$ such that ${\rm in}_{<}(B)$ divides $x_2^{c_{2}-d_2}x_3^{d_3}$; in particular ${\rm in}_{<}(B)=x_2^{a}x_3^{b}$ for some $a \leq c_{2}-d_{2}$ and $b \leq d_{3}$. Note that $a \neq 0$. If $a = 0$, then ${\rm in}_{<}(B)= x_3^b$ with $b < c_3$, contradicting the minimality of $c_3$. Similarly, $b \neq 0$. Moreover, $d_1 \neq 0$ since $d_2 < c_2$, and $d_4 \neq 0$ since $d_3<c_3$. Since $B \in I({\bf a})$ and $I({\bf a})$ is generated by the binomials $f_1$, $f_2$, and $f_3$, the monomial $x_2^{a}x_3^{b}$ must be divisible by one of the monomials $x_2^{c_2}$, $x_3^{c_3}$, $x_1^{d_1}x_2^{d_2}$, or $x_3^{d_3}x_4^{d_4}$, a contradiction. \hfill $\square$

\begin{rem1} {\rm In each of the cases considered in Theorems \ref{ACM2} and \ref{ACM3}, $I^{h}({\bf a})$ is not a complete intersection.}
\end{rem1}

\begin{thm1} \label{ACM4} Suppose that $d_2=0$ and $d_{3}<c_{3}$. Then $\overline{C({\bf a})}$ is arithmetically Cohen--Macaulay if and only if $d_{1} \geq d_{3}+d_{4}$.
\end{thm1}
{\em \noindent Proof.} ($\Leftarrow$) Let $f_{1}=x_2^{c_2}-x_1^{c_1}, f_{2}=x_3^{c_3}-x_4^{c_4}$, and $f_{3}=x_1^{d_1}-x_3^{d_3}x_4^{d_4}$. Then their initial monomials $x_2^{c_2}$, $x_3^{c_3}$, and $x_1^{d_1}$ are pairwise relatively prime, hence $\{f_i \mid 1 \leq i \leq 3 \}$ is a Gr\"obner basis of $I({\bf a})$ with respect to $<$ by Remark~\ref{VeryBasicRemark}(1). Since none of the initial monomials is divisible by $x_4$, $\overline{C({\bf a})}$ is arithmetically Cohen--Macaulay by \cite[Theorem 2.2]{HerSta}.\\
($\Rightarrow$) Assume that $d_{1}<d_{3}+d_{4}$. Since $\overline{C({\bf a})}$ is arithmetically Cohen--Macaulay and ${\rm in}_{<}(f_{3})=x_3^{d_3}x_4^{d_4}$, \cite[Theorem 2.2]{HerSta} yields a binomial $B \in I({\bf a})$ such that ${\rm in}_{<}(B)$ divides $x_3^{d_{3}}$, contradicting the minimality of $c_{3}$. \hfill $\square$ \par\vspace{0.5\baselineskip}

As an application of Theorem \ref{ACM4}, we obtain the following corollary. 

\begin{cor1} \label{CorCase1CI} Suppose that $d_2=0$ and $d_{3}<c_{3}$. Then the ideal $I^{h}({\bf a})$ is a complete intersection if and only if $d_{1} \geq d_{3}+d_{4}$.
\end{cor1}

{\em \noindent Proof.} ($\Leftarrow$) By the proof of Theorem~\ref{ACM4}, $\{f_i \mid 1\le i\le 3\}$ is a Gr\"obner basis for $I({\bf a})$ with respect to $<$. Hence, by Remark~\ref{VeryBasicRemark}(2), $\{f_i^{h} \mid 1\le i\le 3\}$ is a Gr\"obner basis for $I^{h}({\bf a})$ with respect to $<_0$. Thus $I^{h}({\bf a})$ is a complete intersection.\\ ($\Rightarrow$) If $I^{h}({\bf a})$ is a complete intersection, then $\overline{C({\bf a})}$ is arithmetically Cohen--Macaulay. By Theorem \ref{ACM4}, it follows that $d_{1} \geq d_{3}+d_{4}$. \hfill $\square$ \par\vspace{0.5\baselineskip}

We conclude this section with an example illustrating the complete intersection criterion given in Corollary~\ref{CorCase1CI}.

\begin{ex1} {\rm Let ${\bf a}=(14, 7, 8, 20)$ and ${\bf v}=(7,7,4,10)$. By \cite[Theorem 2.6]{Katsa}, for any $t \in \mathbb{N}$ such that ${\rm gcd}(14+7t,7+7t,8+4t,20+10t)=1$, the ideal $I({\bf a}+t{\bf v})$ is generated by $x_{2}^{2+t}-x_{1}^{1+t}$, $x_{3}^{5}-x_{4}^{2}$, and $x_{1}^{2}-x_{3}x_{4}$. Here, $d_1=2$, $d_2=0$, $d_3=d_4=1$, and $c_3=5$. By Corollary \ref{CorCase1CI}, the ideal $I^{h}({\bf a}+t{\bf v})$ is a complete intersection whenever ${\rm gcd}(14+7t,7+7t,8+4t,20+10t)=1$.}
\end{ex1}

\section{Projective closures in case~2}

This section focuses on Case 2 of Theorem \ref{CasesMinimal}, where $I({\bf a})$ has a minimal generating set $\{g_{1}=x_{2}^{c_2}-x_{1}^{c_1}, g_{2}=x_{3}^{c_3}-x_{1}^{c_1}, g_{3}=x_{4}^{c_4}-x_{1}^{d_1}x_{2}^{d_2}x_{3}^{d_3}\}$, and at least two of $d_1$, $d_2$, and $d_3$ are nonzero.

We begin by studying the case where $a_{i} = {\rm max}\{a_1,\ldots,a_4\}$ for some $1 \leq i \leq 3$. {\em Without loss of generality, assume $a_1={\rm max}\{a_1,\ldots,a_4\}$, and fix $d_{2} \leq c_2$ and $d_{3} \leq c_3$.}

Throughout this section, unless stated otherwise, let $<$ denote the degree
reverse lexicographic order on $K[x_1,\ldots,x_4]$ with
$x_4>x_3>x_2>x_1$.

We first consider the case $d_2=c_{2}$ and $d_{3}=c_3$.

\begin{thm1} \label{Veryfirst} If $d_2=c_{2}$ and $d_{3}=c_3$, then $I^{h}({\bf a})$ is a complete intersection.
\end{thm1}

{\em \noindent Proof.} Using $g_1$ and $g_2$, we obtain another generating set $\{g_1,g_2,f\}$ of $I({\bf a})$, where $f=x_{4}^{c_4}-x_{1}^{d_1+2c_1}$. Then ${\rm in}_{<}(g_1)=x_{2}^{c_2}$, ${\rm in}_{<}(g_2)=x_{3}^{c_3}$, and ${\rm in}_{<}(f)=x_{4}^{c_4}$ are pairwise relatively prime. Hence, by Remark~\ref{VeryBasicRemark}(1), $\{g_1, g_2, f\}$ is a Gr\"obner basis for $I({\bf a})$ with respect to $<$. By Remark \ref{VeryBasicRemark}(2), $\{g_{1}^{h}, g_{2}^{h}, f^{h}\}$ is a Gr\"obner basis for $I^{h}({\bf a})$ with respect to $<_{0}$, so $I^{h}({\bf a})$ is a complete intersection. \hfill $\square$ \par\vspace{0.5\baselineskip}

To illustrate Theorem \ref{Veryfirst}, we present the following example.

\begin{ex1} {\rm Let ${\bf a}=(6851, 5797, 4433, 4862)$. Then $I({\bf a})$ is generated by $x_{2}^{13}-x_1^{11}$, $x_{3}^{17}-x_{1}^{11}$, and $x_{4}^{31}-x_{2}^{13}x_{3}^{17}$. Here, $c_2=d_2=13$ and $c_3=d_3=17$. By Theorem \ref{Veryfirst}, $I^{h}({\bf a})$ is a complete intersection.}
\end{ex1}

Next, we examine the case $d_2<c_2$ and $d_3=c_3$.

\begin{thm1} \label{ACM5} If $d_2<c_{2}$ and $d_{3}=c_3$, then $\overline{C({\bf a})}$ is arithmetically Cohen--Macaulay if and only if $c_{4} \geq c_1+d_{1}+d_{2}$.
\end{thm1}

{\em \noindent Proof.} ($\Leftarrow$) Let $f=x_{4}^{c_4}-x_{1}^{d_1+c_1}x_{2}^{d_2} \in I({\bf a})$. Since $g_3=f-x_1^{d_1}x_2^{d_2}g_2$, $I({\bf a})$ is generated by the binomials $g_1$, $g_2$, and $f$. From $c_{4} \geq c_1+d_{1}+d_{2}$, we have ${\rm in}_{<}(f)=x_{4}^{c_4}$. Then ${\rm in}_{<}(g_1)$, ${\rm in}_{<}(g_2)$, and ${\rm in}_{<}(f)$ are pairwise relatively prime, so $\{g_1, g_2, f\}$ is a Gr\"obner basis for $I({\bf a})$ with respect to $<$ by Remark~\ref{VeryBasicRemark}(1). As $x_1$ divides none of these initial monomials, $\overline{C({\bf a})}$ is arithmetically Cohen--Macaulay by \cite[Theorem 2.2]{HerSta}.\\
($\Rightarrow$) Assume that $c_{4}<c_1+d_{1}+d_{2}$. Since $\overline{C({\bf a})}$ is arithmetically Cohen--Macaulay and ${\rm in}_{<}(f)=x_{1}^{d_1+c_1}x_{2}^{d_2}$, \cite[Theorem~2.2]{HerSta} guarantees the existence of a binomial $B \in I({\bf a})$ such that ${\rm in}_{<}(B)$ divides $x_{2}^{d_2}$. Thus ${\rm in}_{<}(B)=x_{2}^{e}$ for some $e \leq d_2$, which implies $e<c_2$, contradicting the minimality of $c_2$. \hfill $\square$ \par\vspace{0.5\baselineskip}

The following corollary is a consequence of Theorem \ref{ACM5}.

\begin{cor1} \label{Cor12} If $d_2<c_{2}$ and $d_{3}=c_3$, then $I^{h}({\bf a})$ is a complete intersection if and only if $c_{4} \geq c_1+d_{1}+d_{2}$.
\end{cor1}

{\em \noindent Proof.} If $I^{h}({\bf a})$ is a complete intersection, then $\overline{C({\bf a})}$ is arithmetically Cohen--Macaulay, and hence $c_{4} \geq c_1+d_{1}+d_{2}$ by Theorem \ref{ACM5}. Conversely, assume that $c_{4} \geq c_1+d_{1}+d_{2}$. By the proof of Theorem~\ref{ACM5}, $\{g_1,g_2,f\}$ is a Gr\"obner basis for $I({\bf a})$ with respect to $<$. Hence, by Remark~\ref{VeryBasicRemark}(2), $\{g_1^{h},g_2^{h},f^{h}\}$ is a Gr\"obner basis for $I^{h}({\bf a})$ with respect to $<_0$. Consequently, $I^{h}({\bf a})$ is a complete intersection. \hfill $\square$ \par\vspace{0.5\baselineskip}

The following example illustrates Corollary \ref{Cor12}.

\begin{ex1} {\rm Let ${\bf a}=(60, 20, 15, 16)$. Then $I({\bf a})$ is generated by $x_{2}^3-x_1$, $x_{3}^{4}-x_{1}$, and $x_{4}^{5}-x_{2}x_{3}^4$. Here, $d_2=1<3=c_2$ and $d_3=c_3=4$. Since $c_4=5>2=c_1+d_1+d_2$, Corollary~\ref{Cor12} implies that
$I^{h}({\bf a})$ is a complete intersection.}
\end{ex1}

We now consider the case $d_2=c_2$ and $d_3<c_3$.

\begin{thm1} \label{ACM6} If $d_2=c_{2}$ and $d_{3}<c_3$, then $\overline{C({\bf a})}$ is arithmetically Cohen--Macaulay if and only if $c_{4} \geq d_1+c_1+d_{3}$.
\end{thm1}

{\em \noindent Proof.} ($\Leftarrow$) Let $f=x_{4}^{c_4}-x_{1}^{d_1+c_1}x_{3}^{d_3}\in I({\bf a})$. Since $g_3=f-x_1^{d_1}x_3^{d_3}g_1$, $I({\bf a})$ is generated by $g_1$, $g_2$, and $f$.
From $c_{4} \geq d_1+c_1+d_{3}$, we have ${\rm in}_{<}(f)=x_{4}^{c_4}$.
Then ${\rm in}_{<}(g_1)$, ${\rm in}_{<}(g_2)$, and ${\rm in}_{<}(f)$ are pairwise relatively prime, so $\{g_1, g_2, f\}$ is a Gr\"obner basis for $I({\bf a})$ with respect to $<$ by Remark~\ref{VeryBasicRemark}(1). As $x_1$ divides none of these initial monomials, $\overline{C({\bf a})}$ is arithmetically Cohen--Macaulay by \cite[Theorem 2.2]{HerSta}.\\
($\Rightarrow$) Assume that $c_{4}<c_1+d_{1}+d_{3}$. Since $\overline{C({\bf a})}$ is arithmetically Cohen--Macaulay and ${\rm in}_{<}(f)=x_{1}^{d_1+c_1}x_{3}^{d_3}$, \cite[Theorem 2.2]{HerSta} guarantees the existence of a binomial $B \in I({\bf a})$ such that ${\rm in}_{<}(B)$ divides $x_{3}^{d_3}$. Thus ${\rm in}_{<}(B)=x_{3}^{e}$ for some $e \leq d_3$, which implies $e<c_3$, contradicting the minimality of $c_3$. \hfill $\square$ \par\vspace{0.5\baselineskip}

Theorem \ref{ACM6} yields the following complete intersection criterion.

\begin{cor1} \label{Cor13} If $d_2=c_{2}$ and $d_{3}<c_3$, then $I^{h}({\bf a})$ is a complete intersection if and only if $c_{4} \geq c_1+d_{1}+d_{3}$.
\end{cor1}

{\em \noindent Proof.} If $I^{h}({\bf a})$ is a complete intersection, then $\overline{C({\bf a})}$ is arithmetically Cohen--Macaulay. By Theorem~\ref{ACM6}, it follows that $c_4 \ge c_1+d_1+d_3$. Conversely, assume that $c_{4} \geq c_1+d_{1}+d_{3}$. By the proof of Theorem~\ref{ACM6}, $\{g_1,g_2,f\}$ is a Gr\"obner basis for $I({\bf a})$ with respect to $<$. Hence, by Remark~\ref{VeryBasicRemark}(2), $\{g_1^{h},g_2^{h},f^{h}\}$ is a Gr\"obner basis for $I^{h}({\bf a})$ with respect to $<_0$. Thus, $I^{h}({\bf a})$ is a complete intersection. \hfill $\square$ \par\vspace{0.5\baselineskip}

We now study the case $d_2<c_2$ and $d_3<c_3$.

\begin{thm1} \label{Basic4} Suppose that $d_2<c_2$ and $d_{3}<c_3$. Then the ideal $I^{h}({\bf a})$ is a complete intersection if and only if $c_{4} \geq d_{1}+d_{2}+d_{3}$.
\end{thm1}

{\em \noindent Proof.} ($\Leftarrow$) Since $c_{4} \geq d_1+d_2+d_{3}$, we have ${\rm in}_{<}(g_3)=x_{4}^{c_4}$. The initial monomials of $g_1$, $g_2$, and $g_3$ are therefore pairwise
relatively prime. Hence, by Remark~\ref{VeryBasicRemark}(1),
$\{g_i \mid 1 \leq i \leq 3\}$ is a Gr\"obner basis for
$I({\bf a})$ with respect to $<$.
By Remark~\ref{VeryBasicRemark}(2),
$\{g_i^{h} \mid 1 \leq i \leq 3\}$ is a Gr\"obner basis for
$I^{h}({\bf a})$ with respect to $<_0$.
Therefore $I^{h}({\bf a})$ is a complete intersection.\\ 
($\Rightarrow$) We may choose a minimal generating set $T=\{g_{1}^{h}, g_{2}^{h}, x_{4}^{u_4}-x_{0}^{u_0}x_{1}^{u_1}x_{2}^{u_2}x_{3}^{u_3}\}$ of $I^{h}({\bf a})$, where $u_{4} \geq u_{1}+u_{2}+u_{3}$. Note that $u_{4} \geq c_{4}$ since $x_{4}^{u_4}-x_{1}^{u_1}x_{2}^{u_2}x_{3}^{u_3} \in I({\bf a})$. Assume that $c_{4}<d_1+d_{2}+d_{3}$. Since $g_{3}^{h}=x_{0}^{d_1+d_{2}+d_{3}-c_{4}}x_{4}^{c_4}-x_{1}^{d_1}x_{2}^{d_{2}}x_{3}^{d_{3}} \in I^{h}({\bf a})$ and the set $T$ generates $I^{h}({\bf a})$, the monomial $x_{4}^{u_4}$ divides $x_{4}^{c_4}$, so $u_4 \leq c_4$, and hence $u_4=c_4$. Moreover, since $d_2<c_2$ and $d_3<c_3$, $x_{1}^{u_1}x_{2}^{u_2}x_{3}^{u_3}$ divides $x_{1}^{d_1}x_{2}^{d_{2}}x_{3}^{d_{3}}$, so $u_i \leq d_i$ for each $1 \leq i \leq 3$. Consequently, $x_{4}^{c_4} (x_1^{d_1-u_1}x_2^{d_2-u_2}x_3^{d_3-u_3}-1) \in I({\bf a})$, which implies that $x_1^{d_1-u_1}x_2^{d_2-u_2}x_3^{d_3-u_3}-1 \in I({\bf a})$, since $I({\bf a})$ is prime and does not contain the variable $x_4$. However, $I({\bf a})$ contains no binomial of the form $x_{1}^{v_1} \cdots x_{4}^{v_4}-1$ with $(v_1,\ldots,v_4) \neq (0,\ldots,0)$. Therefore, $d_i=u_i$ for all $1 \leq i \leq 3$, and hence $c_4 \geq d_1+d_2+d_3$, a contradiction. \hfill $\square$ \par\vspace{0.5\baselineskip}

When $d_{2}<c_2$, $d_3<c_3$, and $c_4<d_1+d_2+d_3$, the ideal $I^{h}({\bf a})$ is not a complete intersection. Nevertheless, the following Gr\"obner basis computations show that $\overline{C({\bf a})}$ may still be arithmetically Cohen--Macaulay.

Let $w$ be the smallest positive integer such that $c_2-wd_{2} \leq 0$ or $c_3-wd_{3} \leq 0$, namely $w={\rm min} \{l \in \mathbb{N} \setminus \{0\} \mid c_{2}-ld_{2} \leq 0 \ \textrm{or} \ c_{3}-ld_{3} \leq 0\}$. Since $d_2<c_2$ and $d_3<c_3$, we have $w \geq 2$. For $0 \leq i \leq w-2$ define $$p_{i}=x_{2}^{c_{2}-(i+1)d_{2}}x_{4}^{(i+1)c_{4}}-x_{1}^{c_1}x_{3}^{(i+1)d_{3}}, \ q_{i}=x_{3}^{c_{3}-(i+1)d_3}x_{4}^{(i+1)c_{4}}-x_{1}^{c_1}x_{2}^{(i+1)d_2}.$$ Additionally, define \[
r =
\begin{cases}
x_4^{wc_4}-x_1^{c_1}x_2^{wd_2}x_3^{wd_3-c_3},
& \text{if } c_2-wd_2>0,\\[4pt]
x_4^{wc_4}-x_1^{c_1}x_2^{wd_2-c_2}x_3^{wd_3},
& \text{if } c_2-wd_2 \leq 0.
\end{cases}
\]

\begin{prop1} \label{CMBasic12} Let $d_1=0$, $d_{2}<c_{2}$, $d_{3}<c_{3}$, $c_3>c_2$, and $c_4<d_{2}+d_{3}$. Suppose that the following conditions hold:
\begin{enumerate} \item $(w-1)(c_{4}-d_{2}-d_{3})+c_{2}-c_{1} \geq 0$.
\item $w(c_{4}-d_{2}-d_{3})+c_{3}-c_{1} \geq 0$ when $c_{2}-wd_{2}> 0$.
\item $w(c_{4}-d_{2}-d_{3})+c_{2}-c_{1} \geq 0$ when $c_{2}-wd_{2} \leq 0$.
\end{enumerate}
Then the set \[
\begin{aligned}
T = {} & \{\, g_{1}=x_{2}^{c_2}-x_{1}^{c_1},\;
          g_{2}=x_{3}^{c_3}-x_{1}^{c_1},\;
          f=x_{2}^{d_2}x_{3}^{d_3}-x_{4}^{c_4} \,\} \cup \{\, p_i \mid 0 \le i \le w-2 \,\} \\
      & 
        \cup \{\, q_i \mid 0 \le i \le w-2 \,\}
        \cup \{\, r \,\}
\end{aligned}
\]
is a Gr\"obner basis for $I({\bf a})$ with respect to $<$.
\end{prop1}

{\em \noindent Proof.} For each $0 \leq i \leq w-2$, the binomials $p_{i}$, $q_{i}$, and $r$ belong to the ideal $I({\bf a})$ by direct computation. Here, ${\rm in}_{<}(g_{1})=x_{2}^{c_{2}}$, ${\rm in}_{<}(g_{2})=x_{3}^{c_{3}}$, and ${\rm in}_{<}(f)=x_{2}^{d_{2}}x_{3}^{d_{3}}$. Since $i+1 \leq w-1$ for all $0 \leq i \leq w-2$, it follows that $(i+1)(c_{4}-d_{2}-d_{3})+c_{2}-c_{1} \geq (w-1)(c_{4}-d_{2}-d_{3})+c_{2}-c_{1}$, as $c_{4}-d_{2}-d_{3}<0$. Hence, $(i+1)(c_{4}-d_{2}-d_{3})+c_{2}-c_{1} \geq 0$, which implies that ${\rm in}_{<}(p_{i})=x_{2}^{c_{2}-(i+1)d_{2}}x_{4}^{(i+1)c_{4}}$. Since $c_3>c_2$, it follows that $(w-1)(c_{4}-d_{2}-d_{3})+c_{3}-c_{1}>(w-1)(c_{4}-d_{2}-d_{3})+c_{2}-c_{1}$. Thus $(w-1)(c_{4}-d_{2}-d_{3})+c_{3}-c_{1} >0$, and therefore ${\rm in}_{<}(q_{i})=x_{3}^{c_{3}-(i+1)d_3}x_{4}^{(i+1)c_{4}}$. Finally, ${\rm in}_{<}(r)=x_{4}^{wc_4}$. Since ${\rm in}_{<}(g_{1})$ and ${\rm in}_{<}(g_{2})$ are relatively prime, we deduce that $S(g_{1},g_{2}) \stackrel{T} {\longrightarrow} 0$. Similarly, we have $S(g_{1},q_{i}) \stackrel{T} {\longrightarrow} 0$ and $S(g_{2},p_i) \stackrel{T} {\longrightarrow} 0$ for all $0 \leq i \leq w-2$, as well as $S(g_{1},r) \stackrel{T} {\longrightarrow} 0$, $S(g_{2},r) \stackrel{T} {\longrightarrow} 0$, and $S(f,r) \stackrel{T} {\longrightarrow} 0$. 

We have that $S(g_{1},f)=p_{0} \stackrel{p_{0}} {\longrightarrow} 0$ and $S(g_{2},f)=q_{0} \stackrel{q_{0}} {\longrightarrow} 0$. We now show that $S(g_{1}, p_{i}) \stackrel{T}{\longrightarrow} 0$ for every $0 \leq i \leq w-2$. Observe that $$S(g_{1},p_{i})=x_{1}^{c_1}(x_{2}^{(i+1)d_2}x_{3}^{(i+1)d_3}-x_{4}^{(i+1)c_4}) \stackrel{f} {\longrightarrow} 0.$$ Similarly, we have $S(g_{2},q_{i})=x_{1}^{c_1}(x_{2}^{(i+1)d_2}x_{3}^{(i+1)d_3}-x_{4}^{(i+1)c_4}) \stackrel{f} {\longrightarrow} 0.$ 

For each $0 \leq i \leq w-3$ we have $S(f,p_{i})=x_{1}^{c_1}x_{3}^{(i+2)d_3}-x_{2}^{c_2-(i+2)d_2}x_{4}^{(i+2)c_4} \stackrel{p_{i+1}} {\longrightarrow} 0$. If $c_{2}-wd_{2} \leq 0$, then $S(f,p_{w-2})=x_{1}^{c_1}x_{2}^{wd_{2}-c_{2}}x_{3}^{wd_3}-x_{4}^{wc_4} \stackrel{r} {\longrightarrow} 0$. On the other hand, if $c_{2}-wd_{2}>0$, then $S(f,p_{w-2})=x_{1}^{c_1}x_{3}^{wd_3}-x_{2}^{c_{2}-wd_2}x_4^{wc_4} \stackrel{r} {\longrightarrow} x_1^{c_1}x_3^{wd_3-c_3}(x_{3}^{c_3}-x_{2}^{c_{2}}) \stackrel{g_1} {\longrightarrow} x_1^{c_1}x_3^{wd_3-c_3}(x_{3}^{c_3}-x_{1}^{c_{1}}) \stackrel{g_2} {\longrightarrow} 0.$

For every $0 \leq i \leq w-3$ we have $S(f,q_{i})=x_{1}^{c_1}x_{2}^{(i+2)d_2}-x_{3}^{c_3-(i+2)d_3}x_{4}^{(i+2)c_4} \stackrel{q_{i+1}} {\longrightarrow} 0$. We now consider the $S$--polynomial $S(f,q_{w-2})$. If $c_{3}-wd_{3}>0$, then \[ S(f,q_{w-2})=x_1^{c_1}x_2^{wd_2}-x_3^{c_3-wd_3}x_4^{wc_4}
\stackrel{r}{\longrightarrow}
x_1^{c_1}x_2^{wd_2-c_2}(x_2^{c_2}-x_3^{c_3})
\stackrel{g_1}{\longrightarrow} \] \[
x_1^{c_1}x_2^{wd_2-c_2}(x_1^{c_1}-x_3^{c_3})
\stackrel{g_2}{\longrightarrow}
0.
\]
Suppose that $c_3-wd_3 \leq 0$. Then $S(f,q_{w-2})=x_1^{c_1}x_2^{wd_2}x_3^{wd_3-c_3}-x_4^{wc_4}$. If $c_2-wd_2>0$, then $S(f,q_{w-2})\stackrel{r}{\longrightarrow}0$. If $c_2-wd_2 \leq 0$, then
\[
S(f,q_{w-2})
\stackrel{r}{\longrightarrow}
x_1^{c_1}x_2^{wd_2-c_2}x_3^{wd_3-c_3}(x_2^{c_2}-x_3^{c_3})
\stackrel{g_1}{\longrightarrow} \] \[
x_1^{c_1}x_2^{wd_2-c_2}x_3^{wd_3-c_3}(x_1^{c_1}-x_3^{c_3})
\stackrel{g_2}{\longrightarrow}
0.\]

Next, we show that $S(p_{i},p_{j}) \stackrel{T}{\longrightarrow} 0$, where $0 \leq i<j \leq w-2$. We compute $S(p_{i},p_{j})=x_{1}^{c_1}x_{3}^{(i+1)d_3}(x_{2}^{(j-i)d_{2}}x_{3}^{(j-i)d_{3}}-x_{4}^{(j-i)c_4}) \stackrel{f}{\longrightarrow} 0$. We now verify that $S(p_i,r) \stackrel{T} {\longrightarrow} 0$ for every $0 \leq i \leq w-2$. Suppose first that $c_2-wd_2>0$. Then $$S(p_i,r)=x_{1}^{c_1}x_{3}^{wd_{3}-c_3}(x_{2}^{c_2+(w-i-1)d_{2}}-x_{3}^{c_{3}-(w-i-1)d_{3}}x_{4}^{(w-i-1)c_{4}}) \stackrel{q_{w-i-2}} {\longrightarrow}$$ $$x_{1}^{c_1}x_{2}^{(w-i-1)d_{2}}x_{3}^{wd_{3}-c_3}(x_{2}^{c_2}-x_{1}^{c_{1}}) \stackrel{g_{1}} {\longrightarrow} 0.$$ Now suppose that $c_2-wd_2 \leq 0$. Then $$S(p_i,r)=x_{1}^{c_1}x_{3}^{(i+1)d_3}(x_{2}^{(w-i-1)d_2}x_{3}^{(w-i-1)d_3}-x_{4}^{(w-i-1)c_4}) \stackrel{f} {\longrightarrow} 0.$$

Next, we show that $S(q_{i},q_{j}) \stackrel{T}{\longrightarrow} 0$, where $0 \leq i<j \leq w-2$. We compute $$S(q_{i},q_{j})=x_{1}^{c_1}x_{2}^{(i+1)d_2}(x_{2}^{(j-i)d_{2}}x_{3}^{(j-i)d_{3}}-x_{4}^{(j-i)c_4}) \stackrel{f}{\longrightarrow} 0.$$ We now verify that $S(q_i,r) \stackrel{T} {\longrightarrow} 0$ for every $0 \leq i \leq w-2$. Suppose first that $c_2-wd_2>0$. Then $S(q_i,r)=x_{1}^{c_1}x_{2}^{(i+1)d_2}(x_{2}^{(w-i-1)d_2}x_{3}^{(w-i-1)d_3}-x_{4}^{(w-i-1)c_4}) \stackrel{f} {\longrightarrow}  0$. Now assume that $c_2-wd_2 \leq 0$, so that 
$r=x_4^{wc_4}-x_1^{c_1}x_2^{wd_2-c_2}x_3^{wd_3}$. Then $$S(q_i,r)=x_{1}^{c_1}x_{2}^{wd_{2}-c_2}(x_{3}^{c_3+(w-i-1)d_{3}}-x_{2}^{c_{2}-(w-i-1)d_{2}}x_{4}^{(w-i-1)c_{4}}) \stackrel{p_{w-i-2}} {\longrightarrow}$$ $$x_{1}^{c_1}x_{2}^{wd_{2}-c_2}x_{3}^{(w-i-1)d_{3}}(x_{3}^{c_3}-x_{1}^{c_{1}}) \stackrel{g_{2}} {\longrightarrow} 0.$$

Finally, we show that $S(p_{i},q_{j}) \stackrel{T}{\longrightarrow} 0$, where $0 \leq i, j \leq w-2$. If $i=j$, then 
$S(p_i,q_i)=x_1^{c_1}(x_2^{c_2}-x_3^{c_3}) \stackrel{g_1} {\longrightarrow} x_1^{c_1}(x_1^{c_1}-x_3^{c_3}) \stackrel{g_2} {\longrightarrow} 0$. Suppose that $i<j$. Then $$S(p_{i},q_{j})=x_{1}^{c_1}(x_{2}^{c_2+(j-i)d_2}-x_{3}^{c_3-(j-i)d_3}x_{4}^{(j-i)c_4}) \stackrel{g_{1}}{\longrightarrow}$$ $$x_{1}^{c_1}(x_{1}^{c_1}x_{2}^{(j-i)d_2}-x_{3}^{c_3-(j-i)d_3}x_{4}^{(j-i)c_4}) \stackrel{q_{j-i-1}}{\longrightarrow} 0.$$ Now suppose that $i>j$. Then $$S(p_{i},q_{j})=x_{1}^{c_1}(x_{2}^{c_2+(j-i)d_2}x_{4}^{(i-j)c_4}-x_{3}^{c_3+(i-j)d_3}) \stackrel{p_{i-j-1}}{\longrightarrow}$$ $$x_{1}^{c_1}x_{3}^{(i-j)d_3}(x_{1}^{c_1}-x_{3}^{c_3}) \stackrel{g_{2}}{\longrightarrow} 0.$$ Thus, $T$ is a Gr\"obner basis for $I({\bf a})$ with respect to $<$. \hfill $\square$ \par\vspace{0.5\baselineskip}

\begin{thm1} \label{CM16} Let $d_1=0$, $d_{2}<c_{2}$, $d_{3}<c_{3}$, $c_3>c_2$, and $c_4<d_{2}+d_{3}$. Suppose that the following conditions hold:
\begin{enumerate} \item $(w-1)(c_{4}-d_{2}-d_{3})+c_{2}-c_{1} \geq 0$.
\item $w(c_{4}-d_{2}-d_{3})+c_{3}-c_{1} \geq 0$ when $c_{2}-wd_{2}> 0$.
\item $w(c_{4}-d_{2}-d_{3})+c_{2}-c_{1} \geq 0$ when $c_{2}-wd_{2} \leq 0$.

\end{enumerate}
Then $\overline{C({\bf a})}$ is arithmetically Cohen--Macaulay.
\end{thm1}

{\em \noindent Proof.} By Proposition \ref{CMBasic12}, the set $T$ is a Gr\"obner basis for $I({\bf a})$ with respect to $<$. Since $x_1$ does not divide the initial monomial of any element in $T$, it follows from \cite[Theorem 2.2]{HerSta} that $\overline{C({\bf a})}$ is arithmetically Cohen--Macaulay. \hfill $\square$

\begin{prop1} \label{ACM7}
Let $n \geq 4$ be an integer such that $3$ does not divide $n$, and let
$a_{1}=2n(2n-3)$, $a_{2}=2(2n-3)$, $a_{3}=2n$, and $a_{4}=4n-3$. Then
$\overline{C({\bf a})}$ is arithmetically Cohen--Macaulay, and the minimal
number of generators of $I^{h}({\bf a})$ is equal to $2n$.
\end{prop1}

{\em \noindent Proof.}
Since $3$ does not divide $n$, we have
\[
\gcd(a_{3},a_{4})
=\gcd(2n,4n-3)
=\gcd(2n,3)
=1.
\]
Therefore,
$
\gcd(a_{1},a_{2},a_{3},a_{4})=1.$
For every $n \ge4$, we also have
\[
a_{1}=na_{2}=(2n-3)a_3>a_{4},
\]
so $a_{1}=\max\{a_1,\ldots,a_4\}$.

By \cite[Proposition~3.1]{Katsa}, 
$I({\bf a})$ is a complete intersection ideal generated by
$x_{2}^{n}-x_{1}$, $x_{3}^{2n-3}-x_{1}$, and
$x_{4}^{2}-x_{2}x_{3}^{2}$.
Thus, in the notation of Theorem~\ref{CM16},
$
c_{1}=1$, $c_{2}=n$, $c_{3}=2n-3$, $c_{4}=2$, $d_{1}=0$, $d_{2}=1$, and $d_{3}=2$.
As $n\ge4$, we have
$
d_{2}<c_{2}$, $d_{3}<c_{3}$,
$c_{3}>c_{2}$, and $c_{4}<d_{2}+d_{3}$.
Moreover, $w=n-1$, so $
c_{2}-wd_{2}=1>0$. Furthermore, $
(w-1)(c_{4}-d_{2}-d_{3})+c_{2}-c_{1}=1>0$ and $
w(c_{4}-d_{2}-d_{3})+c_{3}-c_{1}=n-3>0.$
Hence, all the hypotheses of Theorem~\ref{CM16} are satisfied, and
$\overline{C({\bf a})}$ is arithmetically Cohen--Macaulay.

By Proposition~\ref{CMBasic12}, the set
\[
\begin{aligned}
T={}&
\{x_{2}^{n}-x_{1},\;
x_{3}^{2n-3}-x_{1},\;
x_{2}x_{3}^{2}-x_{4}^{2}\} \cup
\{x_{2}^{n-i-1}x_{4}^{2i+2}
-x_{1}x_{3}^{2i+2}\mid 0\leq i\leq n-3\}\\
&\cup
\{x_{3}^{2n-5-2i}x_{4}^{2i+2}
-x_{1}x_{2}^{i+1}\mid 0\leq i\leq n-3\} \cup
\{x_{4}^{2n-2}-x_{1}x_{2}^{n-1}x_{3}\}
\end{aligned}
\]
is a Gr\"obner basis of $I({\bf a})$ with respect to $<$. By
Remark~\ref{VeryBasicRemark}(2), its homogenization
\[
\begin{aligned}
&\{x_2^n-x_0^{n-1}x_1,\;
x_3^{2n-3}-x_0^{2n-4}x_1,\;
x_2x_3^2-x_0x_4^2\}
\cup
\{x_2^{n-i-1}x_4^{2i+2}
-x_0^{n-i-2}x_1x_3^{2i+2}
\mid\\
&0\le i\le n-3\}
\cup
\{x_3^{2n-5-2i}x_4^{2i+2}
-x_0^{2n-5-i}x_1x_2^{i+1}
\mid 0\le i\le n-3\}
\cup
\{x_4^{2n-2}-x_0^{n-3}x_1x_2^{n-1}x_3\}
\end{aligned}
\] is a Gr\"obner
basis of $I^{h}({\bf a})$ with respect to $<_{0}$.
A direct computation shows that this set is a minimal generating set
of $I^{h}({\bf a})$ consisting of $2n$ binomials. \hfill $\square$ \par\vspace{0.5\baselineskip}

{\em Finally, we consider the case $a_{4}={\rm max}\{a_1,\ldots,a_4\}$. Without loss of generality, we fix the assumptions $c_{3}>c_{1}$, $c_{2}>c_{1}$, $d_{2} \leq c_2$, and $d_{3} \leq c_3$.} Under these assumptions, we give a characterization for when $I^{h}({\bf a})$ is a complete intersection.

\begin{thm1} \label{Comp2} The ideal $I^{h}({\bf a})$ is a complete intersection if and only if the element $c_{4}a_4 \in \mathbb{N}\{a_1,\ldots,a_4\}$ admits a factorization $c_{4}a_4=u_{1}a_{1}$ for some  $u_{1} \in \mathbb{N}$ with $u_1>c_1$.
\end{thm1}

{\em \noindent Proof.}
($\Leftarrow$) Suppose that $c_4a_4=u_{1}a_1$ for some positive integer $u_{1}$, and let
$f=x_1^{u_{1}}-x_4^{c_4} \in I({\bf a})$. Since $a_4>a_1$, we have $u_{1}>c_4$. Let $J$ be the ideal generated by
$g_1^h$, $g_2^h$, and $f^h=x_1^{u_{1}}-x_{0}^{u_{1}-c_4}x_4^{c_4}$. Clearly $J \subseteq I^h({\bf a})$. We will show that $I^h({\bf a}) \subseteq J$.

Since $c_1a_1=c_2a_2=c_3a_3$, from
$u_{1}a_1=d_1a_1+d_2a_2+d_3a_3$ we obtain
\[
u_{1}=d_1+\frac{c_1d_2}{c_2}+\frac{c_1d_3}{c_3}.
\]
So
$(u_{1}-d_1)c_2c_3=c_1(d_2c_3+d_3c_2)$. Using $c_{i}a_{i}=c_{j}a_{j}$ for $1 \leq i<j \leq 3$, we have $\gcd(c_i,c_j)=1$ by the minimality of $c_i$. Thus $c_2$ divides $c_1d_2c_3$ and $c_3$ divides $c_1d_3c_2$, and therefore $c_2$ divides $d_2$ and $c_3$ divides $d_3$. Let $d_2=k_2c_2$ and $d_3=k_3c_3$ for some nonnegative integers $k_2$ and $k_3$. Note that at least one of $k_2$ and $k_3$ is positive, since at least two of $d_1$, $d_2$, and $d_3$ are nonzero. Then
$u_{1}=d_1+(k_2+k_3)c_1$, and hence
\[
d_1+d_2+d_3-u_{1}
=k_2(c_2-c_1)+k_3(c_3-c_1)>0,
\]
because $c_2>c_1$ and $c_3>c_1$.

We have $B_{1}=x_2^{k_{2}c_2}
-x_0^{k_2(c_2-c_1)}x_1^{k_2c_1}
\in J$ and $B_{2}=x_3^{k_{3}c_3}
-x_0^{k_3(c_3-c_1)}x_1^{k_3c_1}
\in J$. So
\[x_1^{d_1}x_2^{d_2}x_3^{d_3}
-x_0^{d_1+d_2+d_3-u_1}x_1^{u_1}=
x_1^{d_1}x_3^{d_3}
B_{1}+
x_0^{k_2(c_2-c_1)}x_1^{d_1+k_2c_1}
B_{2}
\in J.
\]
Also,
$x_0^{d_1+d_2+d_3-u_{1}}f^h
=x_0^{d_1+d_2+d_3-u_{1}}x_1^{u_{1}}
-x_0^{d_1+d_2+d_3-c_4}x_4^{c_4}$. Hence,
$g_3^h
=-\bigl(x_1^{d_1}x_2^{d_2}x_3^{d_3}
-x_0^{d_1+d_2+d_3-u_{1}}x_1^{u_{1}}\bigr)-x_0^{d_1+d_2+d_3-u_{1}}f^h$, so $g_3^h \in J$.

Consider the degree reverse lexicographic order $<$ on $K[x_0,\ldots,x_4]$ with
$x_1>x_2>x_3>x_4>x_0$. Then
${\rm in}_{<}(g_1^h)=x_2^{c_2}$,
${\rm in}_{<}(g_2^h)=x_3^{c_3}$, and
${\rm in}_{<}(f^h)=x_1^{u_{1}}$.
These monomials are pairwise relatively prime, so
$\{g_1^h,g_2^h,f^h\}$ is a Gr\"obner basis of $J$. Moreover, none of their initial monomials is divisible by $x_0$. By \cite[Lemma~12.1]{Sturmfels95}, $J:x_0^\infty=J$.

Since $I({\bf a})$ is generated by $g_1$, $g_2$, and $g_3$, it follows from
\cite[Proposition~2.7.7]{Simis} that $I^h({\bf a})$ is the saturation with respect to $x_0$ of the ideal generated by $g_1^h$, $g_2^h$, and $g_3^h$. These binomials belong to $J$. Moreover, $J:x_0^\infty=J$. Hence, $I^h({\bf a})\subseteq J$, and therefore $I^h({\bf a})=J$ is a complete intersection.

($\Rightarrow$) Suppose that $I^h({\bf a})$ is a complete intersection, and let $T=\{g_{1}^{h}, g_{2}^{h}, x_{1}^{u_1}-x_{0}^{u_{0}}x_{2}^{u_2}x_{3}^{u_3}x_{4}^{u_4}\}$ be a minimal generating set of $I^{h}({\bf a})$, where $u_{1} \geq c_{1}$ and $u_{1} \geq u_2+u_3+u_4$. Since $a_{4}={\rm max}\{a_1,\ldots,a_4\}$, we have $c_{4}<d_{1}+d_{2}+d_{3}$, and therefore $x_{0}^{d_{1}+d_{2}+d_{3}-c_{4}}x_{4}^{c_{4}}-x_{1}^{d_{1}}x_{2}^{d_{2}}x_{3}^{d_{3}} \in I^{h}({\bf a})$. The fact that $T$ generates $I^{h}({\bf a})$ implies that the monomial $x_{0}^{u_{0}}x_{2}^{u_2}x_{3}^{u_3}x_{4}^{u_4}$ divides $x_{0}^{d_{1}+d_{2}+d_{3}-c_{4}}x_{4}^{c_{4}}$. Thus, $u_{2}=u_{3}=0$ and $u_{4} \leq c_{4}$.
Since
$x_1^{u_1}-x_4^{u_4}\in I({\bf a})$, we have $u_4\geq c_4$, and therefore $u_4=c_4$. Thus
$c_4a_4=u_1a_1$ and $u_1>c_1$ because $c_1a_1 \neq c_4a_4$.
\hfill $\square$ \par\vspace{0.5\baselineskip}

The next example illustrates Theorem \ref{Comp2}.

\begin{ex1} {\rm Let ${\bf a}=(140, 84, 60, 245)$. Then $I({\bf a})$ is generated by the binomials $x_{2}^{5}-x_{1}^{3}$, $x_{3}^{7}-x_{1}^3$, and $x_{4}^4-x_{1}x_{2}^{5}x_{3}^{7}$. Here, $c_1=3$, $c_2=5$, $c_3=7$, $c_4=4$, $d_1=1$, $d_2=5$, and $d_3=7$. Since $4a_4=7a_1$, Theorem~\ref{Comp2} implies that $I^h({\bf a})$ is a complete intersection.}
\end{ex1}

\section{Projective closures in case 3}

In this section, we study Case 3 of Theorem \ref{CasesMinimal}, where $I({\bf a})$ has a minimal generating set $\{x_{2}^{c_2}-x_{1}^{c_1}, x_{3}^{c_3}-x_{1}^{b_1}x_{2}^{b_2}, x_{4}^{c_4}-x_{1}^{d_1}x_{2}^{d_2}x_{3}^{d_3}\}$. Both $b_1$ and $b_2$ are nonzero, and at least two of $d_1$, $d_2$, and $d_3$ are nonzero.

Suppose first that $a_{i}={\rm max}\{a_1,\ldots,a_4\}$ for some $i \in \{1,2\}$. {\em Without loss of generality, we assume $a_{1}={\rm max}\{a_1,\ldots,a_4\}$ and $b_2 \leq c_2$.} The next theorem provides a complete intersection criterion for $I^{h}({\bf a})$ when $b_2<c_2$.

\begin{thm1} \label{Basic5} Suppose that $b_2<c_2$. Then $I^{h}({\bf a})$ is a complete intersection if and only if $c_{3} \geq b_{1}+b_{2}$ and the element $c_{4}a_{4} \in \mathbb{N}\{a_1,\ldots,a_4\}$ admits a factorization $c_{4}a_{4}=v_{1}a_{1}+v_{2}a_{2}+v_{3}a_{3}$ with $c_{4} \geq v_1+v_{2}+v_{3}$, for some $v_1, v_2, v_3 \in \mathbb{N}$. 
\end{thm1}

{\em \noindent Proof.} ($\Leftarrow$) Let $g_{1}=x_{2}^{c_2}-x_{1}^{c_1}, g_{2}=x_{3}^{c_3}-x_{1}^{b_1}x_{2}^{b_2}$, and $g_{3}=x_{4}^{c_4}-x_{1}^{v_1}x_{2}^{v_2}x_{3}^{v_3}$, where $c_{4} \geq v_1+v_{2}+v_{3}$. Note that $g_3$ is a minimal generator of $I({\bf a})$ by the minimality of $c_4$. Consider the degree reverse lexicographic order $<$ with $x_{4}>x_{3}>x_{2}>x_{1}$. Since $I({\bf a})$ is generated by $g_1$, $g_2$, and $g_3$ and their initial monomials ${\rm in}_{<}(g_1)=x_{2}^{c_2}$, ${\rm in}_{<}(g_2)=x_{3}^{c_3}$, and ${\rm in}_{<}(g_{3})=x_{4}^{c_4}$ are pairwise relatively prime, Remark \ref{VeryBasicRemark}(1) implies that $\{g_{i} \mid 1 \leq i \leq 3\}$ is a Gr\"obner basis for $I({\bf a})$ with respect to $<$. By Remark \ref{VeryBasicRemark}(2), $\{g_{i}^{h} \mid 1 \leq i \leq 3\}$ is a Gr\"obner basis for $I^{h}({\bf a})$ with respect to $<_{0}$. Therefore $I^{h}({\bf a})$ is a complete intersection.\\
($\Rightarrow$) Assume first that $c_3<b_1+b_2$. Since $\overline{C({\bf a})}$ is arithmetically Cohen--Macaulay and ${\rm in}_{<}(g_2)=x_1^{b_1}x_2^{b_2}$, \cite[Theorem~2.2]{HerSta} yields a binomial $B\in I({\bf a})$
whose initial monomial divides $x_2^{b_2}$, contradicting
the minimality of $c_2$. Now suppose that for every factorization $c_{4}a_{4}=v_{1}a_{1}+v_{2}a_{2}+v_{3}a_{3}$ we have $c_{4}<v_1+v_{2}+v_{3}$. We may choose a minimal generating set $T=\{g_{1}^{h}, g_{2}^{h}, x_{4}^{z_4}-x_{0}^{z_{0}}x_{1}^{z_1}x_{2}^{z_2}x_{3}^{z_3}\}$ of $I^{h}({\bf a})$, where $z_{4} \geq z_{1}+z_{2}+z_{3}$. Since $x_{4}^{z_4}-x_{1}^{z_1}x_{2}^{z_2}x_{3}^{z_3} \in I({\bf a})$, we have $z_{4} \geq c_4$. If $z_4=c_4$, then $c_4a_4=z_1a_1+z_2a_2+z_3a_3$ with $z_1+z_2+z_3 \le c_4$, a contradiction. Thus $z_4>c_4$. Moreover, $x_{0}^{d_1+d_2+d_3-c_4}x_{4}^{c_4}-x_{1}^{d_1}x_2^{d_2}x_{3}^{d_3} \in I^{h}({\bf a})$, and since $T$ generates $I^{h}({\bf a})$, we have $x_{4}^{z_4}$ divides $x_{4}^{c_4}$, hence $z_{4} \leq c_4$, a contradiction. \hfill $\square$ \par\vspace{0.5\baselineskip}

The failure of either condition in Theorem \ref{Basic5} prevents $I^{h}({\bf a})$ from being a complete intersection. The following example illustrates this situation.

\begin{ex1} {\rm Let ${\bf a}=(12, 4, 10, 11)$. Then $I({\bf a})$ is generated by $x_{2}^{3}-x_1$, $x_{3}^{2}-x_{1}x_{2}^{2}$, and $x_{4}^{2}-x_{2}^{3}x_{3}$. Here, $b_{2}=2<3=c_{2}$, $c_{3}=2<3=b_{1}+b_{2}$, and $2a_4=a_1+a_3$. By Theorem \ref{Basic5}, $I^{h}({\bf a})$ is not a complete intersection.}
\end{ex1}

Next, we provide a complete intersection criterion for $I^{h}({\bf a})$ when $b_2=c_2$.

\begin{thm1} \label{Basic23} Suppose that $b_2=c_2$. Then $I^{h}({\bf a})$ is a complete intersection if and only if the element $c_{4}a_{4} \in \mathbb{N}\{a_1,\ldots,a_4\}$ admits a factorization $c_{4}a_{4}=v_{1}a_{1}+v_{2}a_{2}+v_{3}a_{3}$ for some $v_1, v_2, v_3 \in \mathbb{N}$ with $c_{4} \geq v_1+v_{2}+v_{3}$. 
\end{thm1}

{\em \noindent Proof.} ($\Leftarrow$) Let $g_{1}=x_{2}^{c_2}-x_{1}^{c_1}, g_{2}=x_{3}^{c_3}-x_{1}^{b_1+c_1}$, and $g_{3}=x_{4}^{c_4}-x_{1}^{v_1}x_{2}^{v_2}x_{3}^{v_3}$, where $c_{4} \geq v_1+v_{2}+v_{3}$. Note that $g_3$ is a minimal generator of $I({\bf a})$ by the minimality of $c_4$. Consider the degree reverse lexicographic order $<$ with $x_{4}>x_{3}>x_{2}>x_{1}$. Since $I({\bf a})$ is generated by $g_1$, $g_2$, and $g_3$ and their initial monomials ${\rm in}_{<}(g_1)=x_{2}^{c_2}$, ${\rm in}_{<}(g_2)=x_{3}^{c_3}$, and ${\rm in}_{<}(g_3)=x_{4}^{c_4}$ are pairwise relatively prime, Remark \ref{VeryBasicRemark}(1) implies that $\{g_{i} \mid 1 \leq i \leq 3\}$ is a Gr\"obner basis for $I({\bf a})$ with respect to $<$. By Remark \ref{VeryBasicRemark}(2), $\{g_{i}^{h} \mid 1 \leq i \leq 3\}$ is a Gr\"obner basis for $I^{h}({\bf a})$ with respect to $<_{0}$. Thus, $I^{h}({\bf a})$ is a complete intersection.\\ ($\Rightarrow$) The proof of this direction is similar to that of the $(\Rightarrow)$ direction of Theorem~\ref{Basic5}, relying on the factorization of $c_4a_4$, and is therefore omitted. \hfill $\square$ \par\vspace{0.5\baselineskip}

The following example illustrates Theorem~\ref{Basic23}.

\begin{ex1} {\rm Let ${\bf a}=(18, 9, 12, 14)$. Then $I({\bf a})$ is generated by $x_{2}^{2}-x_{1}$, $x_{3}^{3}-x_{1}x_{2}^{2}$, and $x_{4}^{3}-x_{2}^2x_{3}^{2}$. Here, $b_{2}=c_{2}=2$ and $3a_4=a_1+2a_3$. By Theorem \ref{Basic23}, $I^{h}({\bf a})$ is a complete intersection.} 
\end{ex1}

{\em Suppose now that $a_{3}={\rm max}\{a_1,\ldots,a_4\}$. Without loss of generality, we fix the assumption $c_{2}>c_{1}$.} The following theorem establishes a complete intersection criterion in the case where $c_4 \geq d_1+d_2+d_3$.

\begin{thm1} \label{basicvfff} Suppose that $c_{4} \geq d_1+d_2+d_3$. Then $I^{h}({\bf a})$ is a complete intersection if and only if the element $c_{3}a_3 \in \mathbb{N}\{a_1,\ldots,a_4\}$ admits a factorization $c_{3}a_3=u_1a_{1}$ for some $u_1 \in \mathbb{N}$ with $u_{1}>c_1$.
\end{thm1}
{\em \noindent Proof.} ($\Leftarrow$) Let $f_{1}=x_{2}^{c_2}-x_{1}^{c_1}, f_{2}=x_{4}^{c_4}-x_{1}^{d_1}x_2^{d_2}x_{3}^{d_3}, f_{3}=x_{1}^{u_1}-x_{3}^{c_3}$. Note that $f_3$ is a minimal generator of $I({\bf a})$ by the minimality of $c_3$. Consider the degree reverse lexicographic order $<$ with $x_4>x_2>x_1>x_3$. Since $I({\bf a})$ is generated by $f_1$, $f_2$, and $f_3$ and their initial monomials ${\rm in}_{<}(f_{1})=x_{2}^{c_2}$, ${\rm in}_{<}(f_{2})=x_{4}^{c_4}$, and ${\rm in}_{<}(f_{3})=x_{1}^{u_1}$ are pairwise relatively prime, Remark~\ref{VeryBasicRemark}(1) implies that $\{f_i \mid 1 \leq i \leq 3\}$ is a Gr\"obner basis for $I({\bf a})$ with respect to $<$. By Remark \ref{VeryBasicRemark}(2), $\{f_{i}^h \mid 1 \leq i \leq 3\}$ is a Gr\"obner basis for $I^{h}({\bf a})$ with respect to $<_{0}$. Thus, $I^{h}({\bf a})$ is a complete intersection.\\
($\Rightarrow$) Let $T=\{x_{1}^{u_1}-x_{0}^{u_0}x_{2}^{u_2}x_{3}^{u_3}x_{4}^{u_4}, x_{2}^{c_2}-x_{0}^{c_2-c_1}x_{1}^{c_1}, x_{4}^{c_4}-x_{0}^{c_4-d_1-d_2-d_3}x_{1}^{d_1}x_{2}^{d_2}x_{3}^{d_3}\}$ be a minimal generating set of $I^{h}({\bf a})$, where $u_{1} \geq u_2+u_3+u_4$. Since $a_{3}={\rm max}\{a_1,\ldots,a_4\}$, we have $c_{3}<b_{1}+b_{2}$, and therefore $x_{0}^{b_{1}+b_{2}-c_{3}}x_{3}^{c_3}-x_{1}^{b_1}x_{2}^{b_2} \in I^{h}({\bf a})$. Furthermore, since $T$ generates $I^{h}({\bf a})$, at least one of $x_{0}^{u_0}x_{2}^{u_2}x_{3}^{u_3}x_{4}^{u_4}$ and $x_{0}^{c_4-d_1-d_2-d_3}x_{1}^{d_1}x_{2}^{d_2}x_{3}^{d_3}$ divides $x_{0}^{b_{1}+b_{2}-c_{3}}x_{3}^{c_3}$. We distinguish two cases.\\
(1) Suppose that $x_{0}^{u_0}x_{2}^{u_2}x_{3}^{u_3}x_{4}^{u_4}$ divides $x_{0}^{b_{1}+b_{2}-c_{3}}x_{3}^{c_3}$. Then $u_{2}=u_{4}=0$ and $u_{3} \leq c_{3}$, so $x_{1}^{u_1}-x_{0}^{u_0}x_{3}^{u_3} \in I^{h}({\bf a})$, implying $x_{1}^{u_1}-x_{3}^{u_3} \in I({\bf a})$. Thus $u_{3} \geq c_{3}$ by the minimality of $c_3$, hence $u_{3}=c_{3}$. Since $c_{1}a_1 \neq c_{3}a_3$, we conclude that $u_1>c_1$.\\
(2) Suppose that $x_{0}^{c_4-d_1-d_2-d_3}x_{1}^{d_1}x_{2}^{d_2}x_{3}^{d_3}$ divides $x_{0}^{b_{1}+b_{2}-c_{3}}x_{3}^{c_3}$. Since the latter monomial contains neither $x_1$ nor $x_2$, it follows that $d_1=d_2=0$, contradicting the assumption that at least two among $d_1$, $d_2$, $d_3$ are strictly positive.  \hfill $\square$

The following example illustrates Theorem~\ref{basicvfff}.

\begin{ex1} {\rm Let ${\bf a}=(4845, 4080, 5491, 4883)$. Then $I({\bf a})$ is generated by $x_{2}^{19}-x_{1}^{16}$, $x_{3}^{15}-x_{1}x_{2}^{19}$, and $x_{4}^{17}-x_{1}^{16}x_{3}$. Here $c_4=17=d_1+d_2+d_3$ and $15a_3=17a_1 $. By Theorem \ref{basicvfff}, $I^{h}({\bf a})$ is a complete intersection.}
\end{ex1}

When $c_4<d_1+d_2+d_3$, the complete intersection property of the projective closure depends on the factorizations of both $c_3a_3$ and $c_4a_4$.

\begin{thm1} \label{Basic24} Suppose that $c_4<d_1+d_2+d_3$. Then $I^{h}({\bf a})$ is a complete intersection if and only if at least one of the following conditions holds: \begin{enumerate} \item[(A)] The element $c_{3}a_3 \in \mathbb{N}\{a_1,\ldots,a_4\}$ admits a factorization $c_{3}a_3=u_1a_{1}$ for some integer $u_{1}>c_1$, and the element $c_{4}a_4 \in \mathbb{N}\{a_1,\ldots,a_4\}$ admits a factorization $c_{4}a_4=v_1a_{1}+v_2a_{2}+v_3a_3$ for some $v_1, v_2, v_3 \in \mathbb{N}$ with $c_4 \geq v_1+v_2+v_3$.
\item[(B)] The element $c_{3}a_3 \in \mathbb{N}\{a_1,\ldots,a_4\}$ admits a factorization $c_{3}a_3=v_4a_{4}$ for some integer $v_{4}>c_4$, and the element $c_{4}a_4 \in \mathbb{N}\{a_1,\ldots,a_4\}$ admits a factorization $c_{4}a_4=u_1a_{1}$ for some $u_1 \in \mathbb{N}$ such that $u_1>c_4$ and $u_1>c_1$.
\end{enumerate}
\end{thm1}
{\em \noindent Proof.} $(\Leftarrow)$ Suppose first that condition (A) holds. Let $f_{1}=x_{2}^{c_2}-x_{1}^{c_1}, f_{2}=x_{1}^{u_1}-x_{3}^{c_3}$, and $f_{3}=x_{4}^{c_4}-x_{1}^{v_1}x_{2}^{v_2}x_{3}^{v_3}$, where $c_{4} \geq v_1+v_{2}+v_{3}$. The binomials $f_2$ and $f_3$ are minimal generators of
$I({\bf a})$ by the minimality of $c_3$ and $c_4$,
respectively. Consider the degree reverse lexicographic order $<$ with $x_4>x_2>x_1>x_3$. Since $I({\bf a})$ is generated by $f_1$, $f_2$, and $f_3$ and their initial monomials ${\rm in}_{<}(f_{1})=x_{2}^{c_2}$, ${\rm in}_{<}(f_{2})=x_{1}^{u_1}$, and ${\rm in}_{<}(f_{3})=x_{4}^{c_4}$ are pairwise relatively prime, Remark~\ref{VeryBasicRemark}(1) implies that $\{f_i \mid 1 \leq i \leq 3\}$ is a Gr\"obner basis for $I({\bf a})$ with respect to $<$. By Remark \ref{VeryBasicRemark}(2), $\{f_1^h,f_2^h,f_3^h\}$ is a Gr\"obner basis for $I^{h}({\bf a})$ with respect to $<_{0}$. Thus, $I^{h}({\bf a})$ is a complete intersection. 

Now suppose that condition (B) holds. Consider the binomials $f_{4}=x_{4}^{v_4}-x_{3}^{c_3}$ and $f_{5}=x_{1}^{u_1}-x_{4}^{c_4}$, where $u_1>c_4$. The binomials $f_4$ and $f_5$ are minimal generators of
$I({\bf a})$ by the minimality of $c_3$ and $c_4$,
respectively. Since $I({\bf a})$ is generated by $f_1$, $f_4$, and $f_5$ and their initial monomials ${\rm in}_{<}(f_{1})$, ${\rm in}_{<}(f_{4})=x_{4}^{v_4}$, and ${\rm in}_{<}(f_{5})=x_{1}^{u_1}$ are pairwise relatively prime, Remark \ref{VeryBasicRemark}(1) implies that $\{f_1, f_4, f_5\}$ is a Gr\"obner basis for $I({\bf a})$ with respect to $<$. By Remark \ref{VeryBasicRemark}(2), $\{f_1^h, f_4^h, f_5^h\}$ is a Gr\"obner basis for $I^{h}({\bf a})$ with respect to $<_{0}$. Consequently, $I^{h}({\bf a})$ is a complete intersection.\\
$(\Rightarrow)$ Let $T=\{x_{1}^{u_1}-x_{0}^{u_0}x_{2}^{u_2}x_{3}^{u_3}x_{4}^{u_4}, x_{2}^{c_2}-x_{0}^{c_2-c_1}x_{1}^{c_1}, x_{4}^{v_4}-x_{0}^{v_0}x_{1}^{v_1}x_{2}^{v_2}x_{3}^{v_3}\}$ be a minimal generating set of $I^{h}({\bf a})$, where $u_{1} \geq u_2+u_3+u_4$ and $v_{4} \geq v_1+v_2+v_3$. Since $c_{3}<b_{1}+b_{2}$, the binomial $x_{0}^{b_{1}+b_{2}-c_{3}}x_{3}^{c_3}-x_{1}^{b_1}x_{2}^{b_2}$ belongs to $I^{h}({\bf a})$. Since $T$ generates $I^{h}({\bf a})$, at least one of $x_{0}^{u_0}x_{2}^{u_2}x_{3}^{u_3}x_{4}^{u_4}$ and $x_{0}^{v_0}x_{1}^{v_1}x_{2}^{v_2}x_{3}^{v_3}$ divides $x_{0}^{b_{1}+b_{2}-c_{3}}x_{3}^{c_3}$. We distinguish two cases.\\
(1) Suppose that $x_{0}^{u_0}x_{2}^{u_2}x_{3}^{u_3}x_{4}^{u_4}$ divides $x_{0}^{b_{1}+b_{2}-c_{3}}x_{3}^{c_3}$. Then $u_{2}=u_{4}=0$ and $u_{3} \leq c_{3}$. Hence, $x_{1}^{u_1}-x_{0}^{u_0}x_{3}^{u_3} \in I^{h}({\bf a})$, and therefore $x_{1}^{u_1}-x_{3}^{u_3} \in I({\bf a})$. Thus $u_{3} \geq c_{3}$ and therefore $u_{3}=c_{3}$. Since $c_1a_1 \neq c_3a_3$, we have $u_1>c_1$. Furthermore, the binomial $x_{0}^{d_1+d_2+d_3-c_4}x_{4}^{c_4}-x_{1}^{d_1}x_2^{d_2}x_3^{d_3}$ belongs to $I^{h}({\bf a})$. Since $T$ generates $I^{h}({\bf a})$, the monomial $x_4^{v_4}$ must divide $x_{0}^{d_1+d_2+d_3-c_4}x_{4}^{c_4}$, and therefore $v_4 \leq c_4$. But $x_{4}^{v_4}-x_{1}^{v_1}x_{2}^{v_2}x_{3}^{v_3} \in I({\bf a})$, so $v_{4} \geq c_4$, and thus $v_4=c_4$. Consequently, the condition (A) holds.\\
(2) Suppose that $x_{0}^{v_0}x_{1}^{v_1}x_{2}^{v_2}x_{3}^{v_3}$ divides $x_{0}^{b_{1}+b_{2}-c_{3}}x_{3}^{c_3}$. Then $v_{1}=v_{2}=0$ and $v_{3} \leq c_{3}$, so $x_{4}^{v_4}-x_{0}^{v_0}x_{3}^{v_3} \in I^{h}({\bf a})$, implying $x_{4}^{v_4}-x_{3}^{v_3} \in I({\bf a})$. Thus $v_{3} \geq c_3$, and therefore $v_{3}=c_{3}$. Since $c_{3}a_{3} \neq c_4a_4$, it follows that $v_4>c_{4}$. Next, consider the binomial $x_{0}^{d_1+d_2+d_3-c_4}x_{4}^{c_4}-x_{1}^{d_1}x_2^{d_2}x_3^{d_3} \in I^{h}({\bf a})$. Since $T$ generates $I^{h}({\bf a})$, the monomial $x_{0}^{u_0}x_{2}^{u_2}x_{3}^{u_3}x_{4}^{u_4}$ must divide $x_{0}^{d_1+d_2+d_3-c_4}x_{4}^{c_4}$, so $u_2=u_3=0$ and $u_4 \leq c_4$. Hence, $x_{1}^{u_1}-x_{0}^{u_0}x_{4}^{u_4} \in I^{h}({\bf a})$ with $u_1 \geq u_4$, implying $x_{1}^{u_1}-x_{4}^{u_4} \in I({\bf a})$ with $u_4 \geq c_4$. Thus $u_4=c_4$. So $u_1 \geq c_4$. Since $a_1 \neq a_4$, it follows that $u_1>c_4$, and because $c_1a_1 \neq c_4a_4$, we also have $u_1>c_1$. \hfill $\square$ \par\vspace{0.5\baselineskip}

The following remark shows that the conditions (A) and (B) of Theorem \ref{Basic24} cannot hold at the same time.

\begin{rem1} {\rm Suppose that both conditions (A) and (B) of Theorem \ref{Basic24} hold. From the proof of Theorem \ref{Basic24}, condition (B) implies that the initial ideal ${\rm in}_{<}(I({\bf a}))$ is generated by $x_2^{c_2}$, $x_4^{v_4}$, and $x_1^{u_1}$, where $v_4>c_4$. By condition (A), there exists a binomial $x_4^{c_4}-x_1^{v_1}x_2^{v_2}x_3^{v_3} \in I({\bf a})$ with $c_4 \ge v_1+v_2+v_3$. 
Thus $x_4^{c_4} \in {\rm in}_{<}(I({\bf a}))$, a contradiction since $v_4>c_4$.}
\end{rem1}

The following examples illustrate some of the cases arising in Theorem~\ref{Basic24}.

\begin{ex1} {\rm Let ${\bf a}=(5083, 4641, 8602, 8372)$. Then $I({\bf a})$ is generated by $x_{2}^{23}-x_{1}^{21}$, $x_{3}^{13}-x_{1}x_{2}^{23}$, and $x_{4}^{17}-x_{1}^{7}x_{2}^{23}$. In this case, we have $c_1=21$, $c_3=13$, and $c_4=17<30=d_1+d_2+d_3$. Note that there are no $v_1, v_2, v_3 \in \mathbb{N}$ such that $17a_{4}=v_1a_{1}+v_{2}a_{2}+v_{3}a_{3}$ with $v_1+v_2+v_3 \leq 17$, and thus condition (A) of Theorem~\ref{Basic24} does not hold. Additionally, although $17a_4=28a_1$, the element $13a_3$ does not admit a factorization $13a_3=v_4a_4$ with $v_4>17$, and therefore condition (B) also fails. By Theorem \ref{Basic24}, $I^{h}({\bf a})$ is not a complete intersection.}
\end{ex1}

\begin{ex1} {\rm Let ${\bf a}=(7429, 6783, 21850, 9775)$. Then $I({\bf a})$ is generated by $x_{2}^{23}-x_{1}^{21}$, $x_{3}^{17}-x_{1}^{29}x_{2}^{23}$, and $x_{4}^{19}-x_{1}^{4}x_{2}^{23}$. In this case, we have $c_1=21$, $c_3=17$, and $c_4=19<27=d_1+d_2+d_3$. Note that $17a_3=38a_4$ and $19a_4=25a_1$. Therefore, condition (B) of Theorem~\ref{Basic24} is satisfied, and it follows that $I^{h}({\bf a})$ is a complete intersection.}
\end{ex1}

{\em Suppose now that $a_{4}={\rm max}\{a_1,\ldots,a_4\}$. Without loss of generality, we fix $c_{2}>c_{1}$.} The next theorem gives a complete intersection criterion for $c_{3} \geq b_1+b_2$.

\begin{thm1} \label{Basic25} Suppose that $c_{3} \geq b_1+b_2$. Then $I^{h}({\bf a})$ is a complete intersection if and only if the element $c_{4}a_4 \in \mathbb{N}\{a_1,\ldots,a_4\}$ admits a factorization $c_{4}a_4=u_1a_{1}$ for some $u_1 \in \mathbb{N}$ with $u_{1}>c_1$.
\end{thm1}

{\em \noindent Proof.} ($\Leftarrow$) Let $f_{1}=x_{2}^{c_2}-x_{1}^{c_1}, f_{2}=x_{3}^{c_3}-x_{1}^{b_1}x_2^{b_2}$, and $f_{3}=x_{1}^{u_1}-x_{4}^{c_4} \in I({\bf a})$. Note that $f_3$ is a minimal generator of
$I({\bf a})$ by the minimality of $c_4$. Consider the degree reverse lexicographic order $<$ with $x_4>x_3>x_2>x_1$. Since $I({\bf a})$ is generated by $f_1$, $f_2$, and $f_3$ and their initial monomials ${\rm in}_{<}(f_{1})=x_{2}^{c_2}$, ${\rm in}_{<}(f_{2})=x_{3}^{c_3}$, and ${\rm in}_{<}(f_{3})=x_{1}^{u_1}$ are pairwise relatively prime, Remark~\ref{VeryBasicRemark}(1) implies that $\{f_i \mid 1 \leq i \leq 3\}$ is a Gr\"obner basis for $I({\bf a})$ with respect to $<$. By Remark \ref{VeryBasicRemark}(2), $\{f_{i}^h \mid 1 \leq i \leq 3\}$ is a Gr\"obner basis for $I^{h}({\bf a})$ with respect to $<_{0}$. Hence, $I^{h}({\bf a})$ is a complete intersection.

($\Rightarrow$) Let $T=\{x_{1}^{u_1}-x_{0}^{u_0}x_{2}^{u_2}x_{3}^{u_3}x_{4}^{u_4}, x_{2}^{c_2}-x_{0}^{c_2-c_1}x_{1}^{c_1}, x_{3}^{c_3}-x_{0}^{c_3-b_1-b_2}x_{1}^{b_1}x_{2}^{b_2}\}$ be a minimal generating set of $I^{h}({\bf a})$, where $u_{1} \geq u_2+u_3+u_4$. Since $a_{4}={\rm max}\{a_1,\ldots,a_4\}$, it follows that $c_{4}<d_{1}+d_{2}+d_3$, so $x_{0}^{d_{1}+d_{2}+d_3-c_{4}}x_{4}^{c_4}-x_{1}^{d_1}x_{2}^{d_2}x_{3}^{d_3} \in I^{h}({\bf a})$. Because $T$ generates $I^{h}({\bf a})$, the monomial $x_{0}^{u_0}x_{2}^{u_2}x_{3}^{u_3}x_{4}^{u_4}$ must divide $x_{0}^{d_{1}+d_{2}+d_3-c_{4}}x_{4}^{c_4}$. Hence, $u_{2}=u_{3}=0$ and $u_{4} \leq c_4$, so $x_{1}^{u_1}-x_{0}^{u_0}x_{4}^{u_4} \in I^{h}({\bf a})$, and therefore $x_{1}^{u_1}-x_{4}^{u_4} \in I({\bf a})$. Thus $u_{4} \geq c_{4}$, and hence $u_{4}=c_{4}$. Since $c_{1}a_1 \neq c_{4}a_4$, we conclude that $u_1>c_1$. \hfill $\square$ \par\vspace{0.5\baselineskip}

The following example illustrates Theorem \ref{Basic25}.

\begin{ex1} {\rm Let ${\bf a}=(455, 385, 120, 1092)$. Then $I({\bf a})$ is generated by $x_{2}^{13}-x_{1}^{11}$, $x_{3}^{7}-x_{1}x_{2}$, and $x_{4}^{5}-x_{1}x_2^{13}$. Here, $c_1=11$, $c_{3}=7> 2=b_1+b_2$, and $c_{4}=5$. Also $5a_4=12a_1$. By Theorem \ref{Basic25}, $I^{h}({\bf a})$ is a complete intersection.}
\end{ex1}

Finally, we provide a complete intersection criterion for $c_3<b_1+b_2$.

\begin{thm1} \label{Basic26} Suppose that $c_3<b_1+b_2$. Then $I^{h}({\bf a})$ is a complete intersection if and only if at least one of the following conditions holds: \begin{enumerate} \item[(A)] The element $c_{4}a_4 \in \mathbb{N}\{a_1,\ldots,a_4\}$ admits a factorization $c_{4}a_4=u_1a_{1}$ for some integer $u_{1}>c_1$, and the element $c_{3}a_3 \in \mathbb{N}\{a_1,\ldots,a_4\}$ admits a factorization $c_{3}a_3=v_1a_{1}+v_2a_{2}+v_4a_4$ for some $v_1, v_2, v_4 \in \mathbb{N}$ with $c_3 \geq v_1+v_2+v_4$.
\item[(B)] The element $c_{4}a_4 \in \mathbb{N}\{a_1,\ldots,a_4\}$ admits a factorization $c_{4}a_4=v_3a_{3}$ for some integer $v_{3}>c_3$, and the element $c_{3}a_3 \in \mathbb{N}\{a_1,\ldots,a_4\}$ admits a factorization $c_{3}a_3=u_1a_{1}$ for some $u_1 \in \mathbb{N}$ such that $u_1>c_3$ and $u_1>c_1$.
\end{enumerate}
\end{thm1}

{\em \noindent Proof.} $(\Leftarrow)$ Suppose first that condition (A) holds. Let $f_{1}=x_{2}^{c_2}-x_{1}^{c_1}, f_{2}=x_{1}^{u_1}-x_{4}^{c_4}$, and $f_{3}=x_{3}^{c_3}-x_{1}^{v_1}x_{2}^{v_2}x_4^{v_4}$, where $c_3 \geq v_1+v_2+v_4$. The binomials $f_2$ and $f_3$ are minimal generators of
$I({\bf a})$ by the minimality of $c_4$ and $c_3$,
respectively. Consider the degree reverse lexicographic order $<$ with $x_3>x_1>x_2>x_4$. Since $I({\bf a})$ is generated by $f_1$, $f_2$, and $f_3$ and their initial monomials ${\rm in}_{<}(f_{1})=x_{2}^{c_2}$, ${\rm in}_{<}(f_{2})=x_{1}^{u_1}$, and ${\rm in}_{<}(f_{3})=x_{3}^{c_3}$ are pairwise relatively prime, Remark~\ref{VeryBasicRemark}(1) implies that $\{f_i \mid 1 \leq i \leq 3\}$ is a Gr\"obner basis for $I({\bf a})$ with respect to $<$. By Remark \ref{VeryBasicRemark}(2), $\{f_{i}^h \mid 1 \leq i \leq 3\}$ is a Gr\"obner basis for $I^{h}({\bf a})$ with respect to $<_{0}$. Thus, $I^{h}({\bf a})$ is a complete intersection. 

Suppose now that condition (B) holds. Consider the binomials $f_{4}=x_{3}^{v_3}-x_{4}^{c_4}$ and $f_{5}=x_{1}^{u_1}-x_{3}^{c_3}$, where $u_1 \geq c_3$. The binomials $f_4$ and $f_5$ are minimal generators of
$I({\bf a})$ by the minimality of $c_4$ and $c_3$,
respectively. Since $I({\bf a})$ is generated by $f_1$, $f_4$, and $f_5$ and their initial monomials ${\rm in}_{<}(f_{1})$, ${\rm in}_{<}(f_{4})=x_{3}^{v_3}$, and ${\rm in}_{<}(f_{5})=x_{1}^{u_1}$ are pairwise relatively prime, Remark~\ref{VeryBasicRemark}(1) implies that $\{f_1, f_4, f_5\}$ is a Gr\"obner basis for $I({\bf a})$ with respect to $<$. By Remark \ref{VeryBasicRemark}(2), $\{f_1^h, f_4^h, f_5^h\}$ is a Gr\"obner basis for $I^{h}({\bf a})$ with respect to $<_{0}$. Therefore, $I^{h}({\bf a})$ is a complete intersection.\\
$(\Rightarrow)$ Let $T=\{x_{1}^{u_1}-x_{0}^{u_0}x_{2}^{u_2}x_{3}^{u_3}x_{4}^{u_4}, x_{2}^{c_2}-x_{0}^{c_2-c_1}x_{1}^{c_1}, x_{3}^{v_3}-x_{0}^{v_0}x_{1}^{v_1}x_{2}^{v_2}x_{4}^{v_4}\}$ be a minimal generating set of $I^{h}({\bf a})$, where $u_{1} \geq u_2+u_3+u_4$ and $v_{3} \geq v_1+v_2+v_4$ with $v_3 \geq c_3$. Since $c_{4}<d_{1}+d_{2}+d_3$, the binomial $x_{0}^{d_{1}+d_{2}+d_3-c_{4}}x_{4}^{c_4}-x_{1}^{d_1}x_{2}^{d_2}x_{3}^{d_3}$ belongs to $I^{h}({\bf a})$. Since $T$ generates $I^{h}({\bf a})$, at least one of $x_{0}^{u_0}x_{2}^{u_2}x_{3}^{u_3}x_{4}^{u_4}$ and $x_{0}^{v_0}x_{1}^{v_1}x_{2}^{v_2}x_{4}^{v_4}$ divides $x_{0}^{d_{1}+d_{2}+d_3-c_{4}}x_{4}^{c_4}$. We distinguish two cases.\\ (1) Suppose that $x_{0}^{u_0}x_{2}^{u_2}x_{3}^{u_3}x_{4}^{u_4}$ divides $x_{0}^{d_{1}+d_{2}+d_3-c_{4}}x_{4}^{c_4}$. Then $u_{2}=u_{3}=0$ and $u_{4} \leq c_{4}$, so $x_{1}^{u_1}-x_{0}^{u_0}x_{4}^{u_4} \in I^{h}({\bf a})$, implying $x_{1}^{u_1}-x_{4}^{u_4} \in I({\bf a})$. Thus $u_{4} \geq c_{4}$, and therefore $u_{4}=c_{4}$. Since $c_1a_1 \neq c_4a_4$, we conclude that $u_1>c_1$. Additionally, $x_{0}^{b_1+b_2-c_3}x_{3}^{c_3}-x_{1}^{b_1}x_2^{b_2} \in I^{h}({\bf a})$. Since $T$ generates $I^{h}({\bf a})$, the monomial $x_3^{v_3}$ divides $x_{0}^{b_1+b_2-c_3}x_{3}^{c_3}$, and therefore $v_3 \leq c_3$. Thus $v_3=c_3$, so condition (A) holds.\\
(2) Suppose that $x_{0}^{v_0}x_{1}^{v_1}x_{2}^{v_2}x_{4}^{v_4}$ divides $x_{0}^{d_{1}+d_{2}+d_3-c_{4}}x_{4}^{c_4}$. Then $v_{1}=v_{2}=0$ and $v_{4} \leq c_{4}$, so $x_{3}^{v_3}-x_{0}^{v_0}x_{4}^{v_4} \in I^{h}({\bf a})$, implying $x_{3}^{v_3}-x_{4}^{v_4} \in I({\bf a})$. Thus $v_{4} \geq c_4$, and therefore $v_{4}=c_{4}$. Since $c_{3}a_{3} \neq c_4a_4$, it follows that $v_3>c_{3}$. Now, consider the binomial $x_{0}^{b_1+b_2-c_3}x_{3}^{c_3}-x_{1}^{b_1}x_2^{b_2} \in I^{h}({\bf a})$. Since $T$ generates $I^{h}({\bf a})$, the monomial $x_{0}^{u_0}x_{2}^{u_2}x_{3}^{u_3}x_{4}^{u_4}$ must divide $x_{0}^{b_1+b_2-c_3}x_{3}^{c_3}$, so $u_2=u_4=0$ and $u_3 \leq c_3$. Hence, $x_1^{u_1}-x_{0}^{u_0}x_{3}^{u_3} \in I^{h}({\bf a})$ with $u_1 \geq u_3$, implying $x_1^{u_1}-x_{3}^{u_3} \in I({\bf a})$ with $u_3 \geq c_3$. Thus $u_3=c_3$. Since $a_1 \neq a_3$, it follows that $u_1>c_3$, and because $c_1a_1 \neq c_3a_3$, we also have $u_1>c_1$. Thus condition (B) holds.  \hfill $\square$

\begin{rem1} {\rm Conditions (A) and (B) of Theorem~\ref{Basic26} cannot hold simultaneously.}
\end{rem1}

We conclude this section with an example illustrating condition {\rm (A)} in Theorem~\ref{Basic26}.

\begin{ex1} {\rm Let ${\bf a}=(3553, 2431, 3069, 5168)$. Then $I({\bf a})$ is generated by $x_{2}^{19}-x_{1}^{13}$, $x_{3}^{17}-x_{1}x_{2}^{20}$, and $x_{4}^{11}-x_{1}^{3}x_{2}^{19}$. Here, $c_3=17<21=b_1+b_2$ and $c_4=11$. Since $16a_1=11a_4$ and $17a_3=14a_1+a_2$, condition (A) of Theorem \ref{Basic26} is satisfied. Thus, $I^{h}({\bf a})$ is a complete intersection.}
\end{ex1}

\noindent{\bf Acknowledgements} The author thanks Dr. Akiyoshi Sannai for his comments on an earlier version of this manuscript.\\

\noindent \textbf{Conflict of interest statement}\\  The author declares no conflicts of interest.\\

\noindent \textbf{Data availability statement}\\Data sharing is not applicable to this article as no datasets were generated or analysed during the current study.

\end{document}